\documentclass[10pt]{article}
\usepackage{graphicx}
\usepackage{caption2}
\usepackage{epsfig}
\usepackage{amssymb}
\usepackage{amsmath}
\usepackage{times}
\usepackage{./mhequ}

\newcommand{\reff}[1]{(\ref{#1})}
\newcommand{\real}{{ {\mathbb R}}}

\renewcommand{\L}{{\mathrm L}}
\renewcommand{\H}{{\mathrm H}}
\renewcommand{\d}{{\mathrm d}}
\newcommand{\D}{{\rm D}}
\newcommand{\E}{{\rm E}}
\renewcommand{\S}{{\cal S}}
\newcommand{\T}{{\cal T}}
\newcommand{\M}{{\rm M}}
\renewcommand{\P}{{\mathbb P}}
\newcommand{\Q}{{\mathbb Q}}

\newcommand{\ed}{{\rm e}}
\newcommand{\erf}{\rm erf}
\newcommand{\sign}{{\rm sign}}

\newcommand{\intR}{
\mbox{\footnotesize${\displaystyle\int_{-\infty}^{\infty}}$}
}

\renewcommand{\vector}[2]{
\Big(
\begin{matrix}   
#1\\ #2
\end{matrix}\Big)
}

\newcommand{\myl}[1]{                   
\hspace{-1mm}
\left(\vbox to #1pt{}\right.
\hspace{-1mm}
}

\newcommand{\my}[2]{
\left#1\vbox to #2pt{}\right.
}

\newcommand{\myr}[1]{
\hspace{-1mm}
\left.\vbox to #1pt{}\right)
\hspace{-1mm}
}

\newcommand{\B}[2]{
B\hspace{-1mm}\my{[}{9}^{#1}_{#2}\my{]}{9}\hspace{-1mm}(t)
}

\newcommand{\Bone}[2]{
B_1\hspace{-1mm}\my{[}{9}^{#1}_{#2}\my{]}{9}\hspace{-1mm}(t)
}

\newenvironment{proof}[1][Proof]{\noindent\textit{#1.} }{\ \rule{0.5em}{0.5em}\par}

\makeatletter
\@addtoreset{equation}{section}
\def\theequation{\thesection.\the\c@equation}
\def\newappendix#1{%
        \let\@oldform\@seccntformat%
        \def\@seccntformat##1{Appendix~\csname the##1\endcsname:~}%
        \section{#1}%
        \let\@seccntformat\@oldform%
        }
\makeatother

\newtheorem{theorem}{Theorem}

\newtheorem{definition}[theorem]{Definition}
\newtheorem{remark}[theorem]{Remark}
\newtheorem{lemma}[theorem]{Lemma}
\newtheorem{proposition}[theorem]{Proposition}

\title{Long tails in the long time asymptotics of quasi-linear
hyperbolic-parabolic systems of conservation laws}
\author{Guillaume van Baalen\footnote{Dept. of Mathematics and
Statistics, Boston University}, Nikola Popovi\'c$^{*}$, and C. Eugene
Wayne$^{*}$}
\begin{document}
\maketitle

\begin{abstract}
\noindent The long-time behaviour of solutions of systems of
conservation laws has been extensively studied. In particular,
Liu and Zeng \cite{liu:1997} have given a detailed exposition of
the leading order asymptotics of solutions close to a constant
background state. In this paper, we extend the analysis of
\cite{liu:1997} by examining higher order terms in the
asymptotics in the framework of the so-called two dimensional
{\em p-system}, though we believe that our methods and results
also apply to more general systems. We give a constructive
procedure for obtaining these terms, and we show that their
structure is determined by the interplay of the parabolic and
hyperbolic parts of the problem. In particular, we prove that the
corresponding solutions develop {\em long tails} that precede the
characteristics.
\end{abstract}

\section{Introduction}

In this paper, we consider the long-time behavior of solutions of
systems of viscous conservation laws. This topic has been
extensively studied. In particular, for the case of solutions
close to a constant background state, \cite{liu:1997} contains a
detailed exposition of the leading order long-time behavior of
such solutions. More precisely, it is shown in \cite{liu:1997}
that the leading order asymptotics are given as a sum of
contributions moving with the characteristic speeds of the
undamped system of conservation laws and that each contribution
evolves as either a Gaussian solution of the heat equation or as
a self-similar solution of the viscous Burger's equation. Thus
with the exception of the translation along characteristics,
these leading order terms reflect primarily the dissipative
aspects of the problem.

In this paper, in an effort to better understand the interplay
between the hyperbolic and parabolic aspects of the problem, we
examine higher order terms in the asymptotics. We work with a
specific two-dimensional system of equations -- the {\em p-system},
but we believe that its behavior is prototypical. In particular,
we think that our methods and results would extend to more
complicated systems such as the `full gas dynamics' and the
equations of Magneto-Hydro-Dynamics (MHD) as considered in
\cite{liu:1997}.

The specific set of equations we consider is the following:
\begin{equa}[2]\label{eqn:p-system}
\partial_t a &= c_1 \partial_x b~, & a(x,0)&=a_0(x)~,\\
\partial_t b &= c_2 \partial_x a + \partial_x g(a,b) + \alpha
\left( \partial_x^2 b + \partial_x( f(a,b) \partial_x b )
\right)~,~~~~& b(x,0)&=b_0(x)~.
\end{equa}

We will make precise the assumptions on the nonlinear terms $f$
and $g$ below, but in order to describe our results informally,
we basically assume that $|g(a,b)| \sim {\cal O}(
(|a|+|b|)^2 )$ and $|f(a,b)| \sim {\cal O}( (|a|+|b|))$. We also
note that without loss of generality, we can set $c_1=c_2=1$ and
$\alpha=2$ in (\ref{eqn:p-system}), which can be achieved by
appropriate scalings of space, time and the dependent variables,
and possible redefinition of the functions $f$ and $g$.

Physically, \reff{eqn:p-system} is a model for compressible,
constant entropy flow, where $a$ represents the volume fraction
(i.e. the reciprocal of the density) and $b$ is the fluid
velocity. The first of the two equations in (\ref{eqn:p-system})
is the consistency relation between these two physical
quantities. In particular, it would not be physically reasonable
to include a dissipative term in this equation, whereas such a
term arises naturally in the second equation which is essentially
Newton's law, in which internal frictional forces are often
present. As a consequence of the form of the dissipation the
damping here is not `diagonalizable' in the terminology of
\cite{liu:1997}.

Next, we note that with the scaling $c_1=c_2=1$ and $\alpha=2$ in
(\ref{eqn:p-system}), the characteristic speeds are $\pm 1$.
Then, following Liu and Zeng \cite{liu:1997}, we introduce new
dependent variables $u$ and $v$ which translate with those
characteristic speeds $\pm1$, respectively. If the initial
conditions $a_0$ and $b_0$ in (\ref{eqn:p-system}) decay
sufficiently fast as $|x|\to\infty$, Liu and Zeng showed that in
the translating frame of reference, $u(x,t)=\frac{1}{\sqrt{1+t}}
g_0(\frac{x}{\sqrt{1+t}}) + {\cal O}((1+t)^{-\frac{3}{4}})$, and
similarly for $v$, where $g_0$ is a self-similar solution of
either the heat equation, or of Burger's equation, depending on
the detailed form of the nonlinear terms. In this paper we derive
similar expressions for the higher order terms in the asymptotics
through a constructive procedure that can be carried out to
arbitrary order.

\begin{figure}[t]
\unitlength=1mm
\begin{center}
\begin{picture}(0,0)(0,0)
\put(-70,-90){\psfig{file=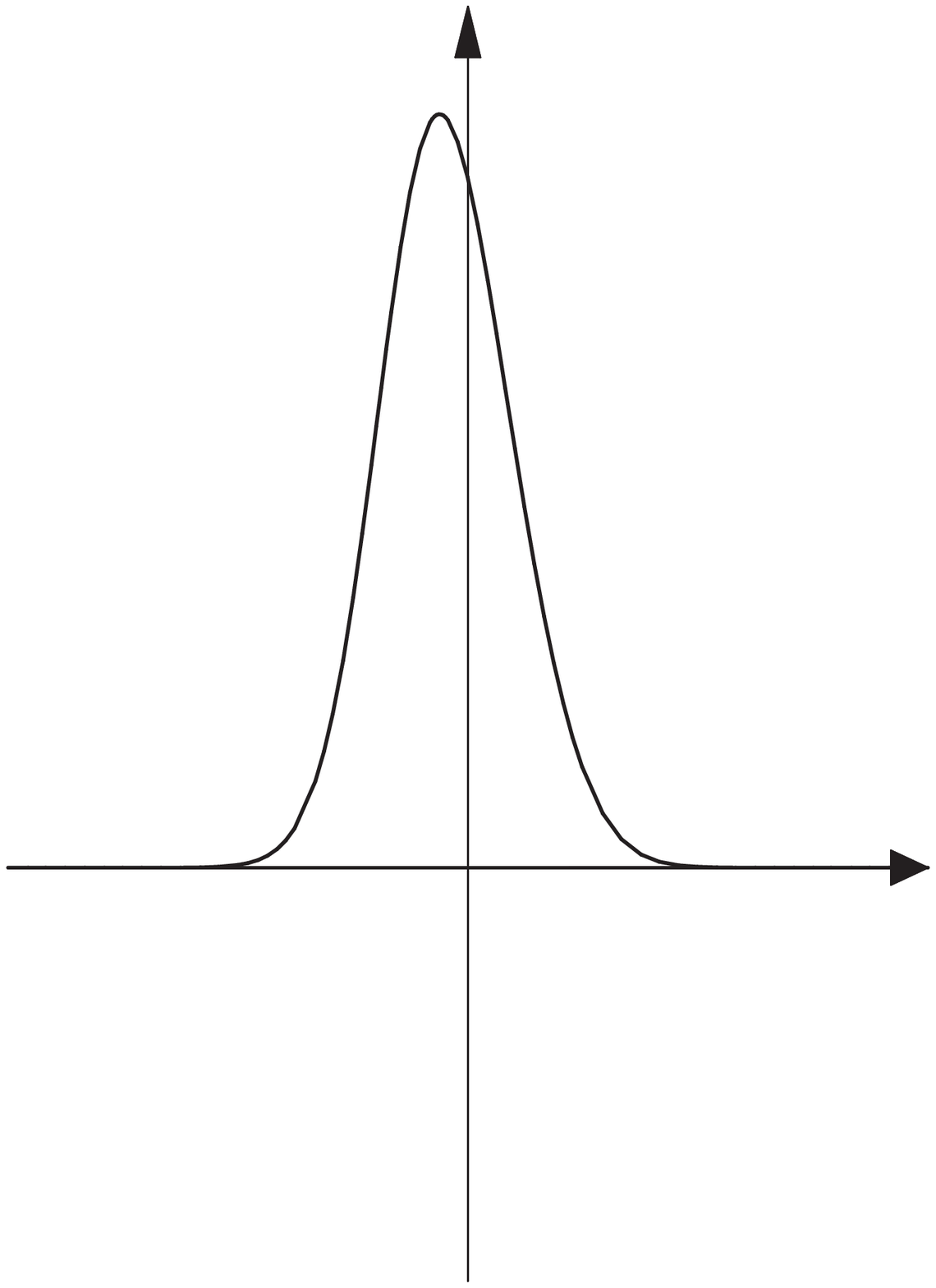,width=6cm}}
\put(0,-90){\psfig{file=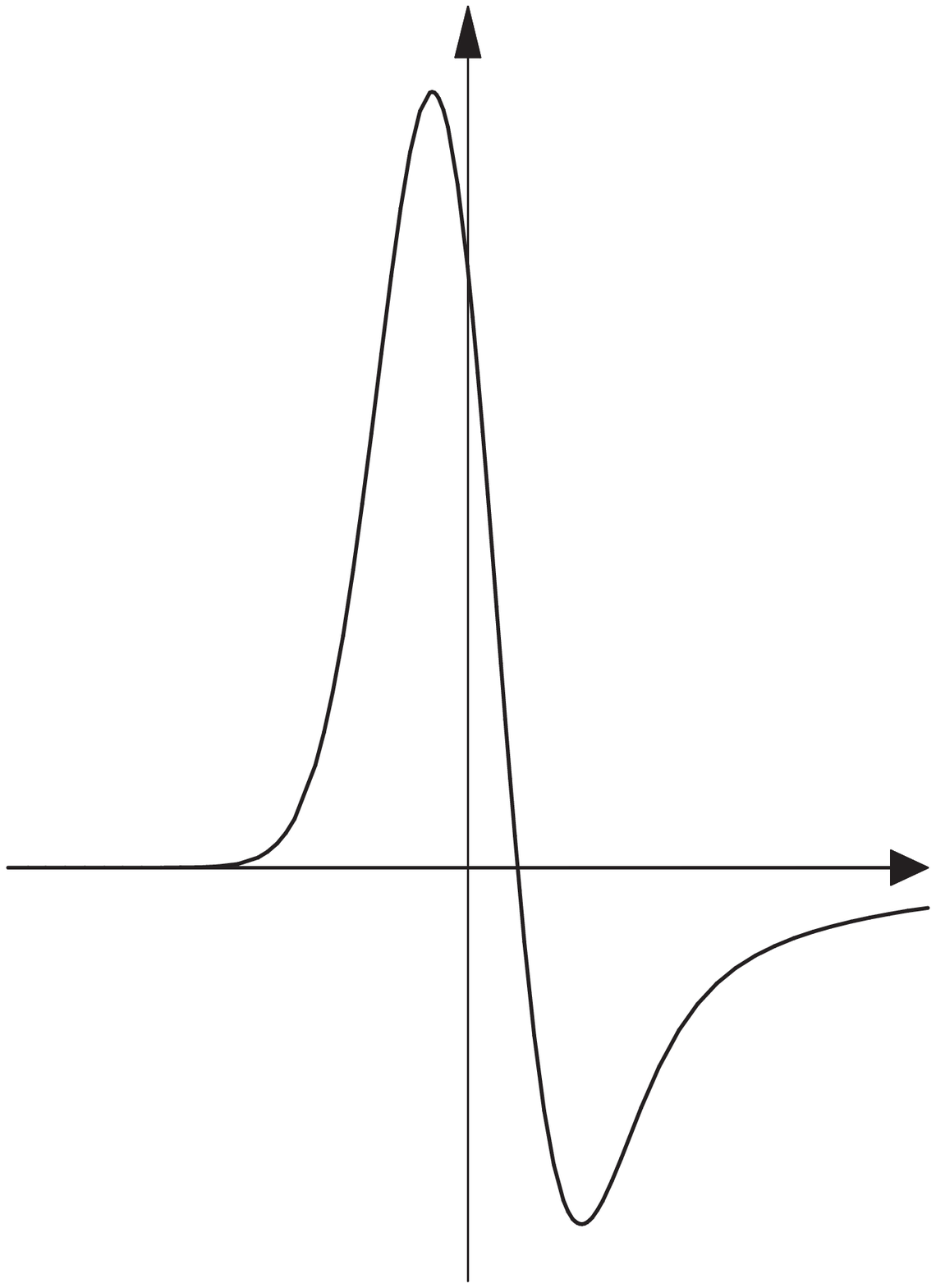,width=6cm}}
\end{picture}
\end{center}
\begin{center}
\begin{picture}(120,80)(0,0)
\put(47,25){$z$}
\put(23,75){$g_0^{+}(z)$}
\put(117,25){$z$}
\put(93,75){$g_1^{+}(z)$}
\end{picture}
\end{center}
\setcaptionwidth{130mm}
\caption{Graphs of the functions $g_0^{+}$ (left panel) and
$g_1^{+}$ (right panel). Note the {\em long tail} of $g_{1}^{+}$
as $z\to\infty$.}
\label{fig:thefigure}
\end{figure}

More precisely, we show that for any $N \ge 1$, there exist
(universal) functions $\{ g_{n}^{\pm}\}_{n=1}^N$ and constants
$\{ d_n^{\pm} \}_{n=1}^N$ determined by the initial conditions,
such that 
\begin{equa}[2]\label{asym}
u(x,t) &= \frac{1}{\sqrt{1+t}}
g_0^{+}({\textstyle\frac{x}{\sqrt{1+t}}})
+ \sum_{n=1}^N 
\frac{1}{(1+t)^{1-\frac{1}{2^{n+1}}}} 
d_n^{+} g_{n}^{+}({\textstyle\frac{x}{\sqrt{1+t}}})
+{\cal O}\myl{12}\frac{1}{(1+t)^{1-\frac{1}{2^{N+2}}}}\myr{12}~.
\end{equa}
We give explicit expressions for the functions $g_{n}^{\pm}$
below, but focusing for the moment on the case $N=1$ and the
variable $u$, we have
\begin{equs}
u(x,t) &= \frac{1}{\sqrt{1+t}}
g_0^{+}({\textstyle\frac{x}{\sqrt{1+t}}})
+ \frac{1}{(1+t)^{\frac{3}{4}}}
d_1^{+} g_{1}^{+} ({\textstyle\frac{x}{\sqrt{1+t}}}) 
+ {\cal O}\myl{12}\frac{1}{(1+t)^{\frac{7}{8}}}\myr{12}~,
\end{equs}
where the functions $g_{0}^{+}(z)$ and $g_{1}^{+}(z)$ are
solutions of the following ordinary differential equations:
\begin{equs}
\partial_z^2 g_{0}^{+}(z) 
+ \frac{1}{2} z \partial_z g_{0}^{+}(z)
+ \frac{1}{2} g_{0}^{+}(z) 
+ c_{+} \partial_z (g_0^{+}(z)^2) &= 0
\label{eqn:nonlineareq}\\
\partial_z^2 g_{1}^{+}(z) 
+ \frac{1}{2} z \partial_z g_{1}^{+}(z)
+ \frac{3}{4} g_{1}^{+}(z) 
+ 2 c_{+} \partial_z (g_0^{+}(z)g_{1}^{+}(z)) &= 0~.
\label{eqn:lineareq}
\end{equs}
Here $c_{+}$ is a constant that depends on the Hessian matrix of
$g(a,b)$ at $a=b=0$ and that will be specified in the course of
our analysis. We will prove that while all solutions of
(\ref{eqn:nonlineareq}) have Gaussian decay as $|x|\to\infty$,
general solutions of the {\em linear} equation
(\ref{eqn:lineareq}) are linear combinations of two functions
$g_{1,\pm}^{+}(z)$, where $g_{1,\pm}^{+}(z)$ decays like a Gaussian
as $z \to \mp \infty$ but only like $|z|^{-\frac{3}{2}}$ as $z \to
\pm \infty$. The graphs of the functions $g_{0}^{+}(z)$ and
$g_{1}^{+}(z)$ are presented in Figure \ref{fig:thefigure}.

Thus, the higher order terms in the asymptotics develop {\em long
tails}. These tails are a manifestation of the hyperbolic part of
the problem (or perhaps more precisely of the interplay between
the parabolic and hyperbolic parts). Were we to consider just the
asymptotic behavior of the viscous Burger's equation which gives
the leading order behavior of the solutions, we would find that
if the initial data is well localized, the higher order terms in
the long-time asymptotics decay rapidly in space and have
temporal decay rates given by half-integers.

Another somewhat surprising aspect of our analysis is that the
tails actually {\em precede} the characteristics.

We also note one additional fact about the expansion in
\reff{asym}. Prior research \cite{gallay:1998,wayne:1997} has
shown that for both parabolic equations and damped wave equations
the eigenfunctions of the operator
\begin{equs}
{\cal L} u(z)  = \partial_z^2 u + \frac{1}{2} z \partial_z u 
\end{equs}
play an important role for the asymptotics. In particular, on
appropriate function
spaces this operator has a sequence of isolated eigenvalues whose
associated eigenfunctions can be used to construct an expansion
for the long-time asymptotics. In this connection we prove that
the functions $g_{n}^{\pm}$ are closely approximated by
eigenfunctions of ${\cal L}$ with eigenvalues $\lambda_n =
-\frac{1}{2} + 2^{-(n+1)}$; more precisely, the functions
$g_{n}^{\pm}$ are eigenfunctions of a compact perturbation of
${\cal L}$, see e.g. (\ref{eqn:lineareq}). However, so far we
have not succeeded in finding a function space which both
contains these eigenfunctions (the functions $g_{n}^{\pm}$ decay
slowly as $z\to\pm\infty$) and in which the corresponding
eigenvalues are isolated points in the spectrum. We plan to
investigate this point further in future research.

Before moving to a precise statement of our results we note that
our approach makes no use of Kawashima's energy estimates for
hyperbolic-parabolic conservation laws \cite{kawashima:1987}. Instead
we prove existence by directly studying the integral form of
\reff{eqn:p-system}.

We now state our results on the Cauchy problem
\reff{eqn:p-system}. We begin by stating the precise assumptions
we make on the nonlinearities $f$ and $g$ in \reff{eqn:p-system}.

\begin{definition}
\label{def:admissible}
The maps $f,g:{\bf R}^2\to{\bf R}$ are admissible nonlinearities
for (\ref{eqn:p-system}) if there is a quadratic map $g_0:{\bf
R}^2\to{\bf R}$ and a constant $C$ such that for all $|{\bf z}|$,
$|{\bf z}_1|$ and $|{\bf z}_2|$ small enough,
\begin{equs}[2]
|g({\bf z})|&\leq C|{\bf z}|^2~,
&~~
|g({\bf z}_1)-g({\bf z}_2)|
&\leq C|{\bf z}_1-{\bf z}_2|(|{\bf z}_1|+|{\bf z}_2|)~,\\
|\Delta g({\bf z})|&\leq C|{\bf z}|^3~,
&~~
|\Delta g({\bf z}_1)-\Delta g({\bf z}_2)|
&\leq C|{\bf z}_1-{\bf z}_2|(|{\bf z}_1|+|{\bf z}_2|)^2~,\\
|f({\bf z})|&\leq C|{\bf z}|~&~\mbox{ and }~~
|f({\bf z}_1)-f({\bf z}_2)|&\leq C|{\bf z}_1-{\bf z}_2|~,
\end{equs}
where $\Delta g({\bf z})\equiv g({\bf z})-g_0({\bf z})$. 
\end{definition}

The main result of this paper can be formulated as follows:
\begin{theorem}\label{thm:maintheorem}
Fix $N>0$. There exists $\epsilon_0>0$ sufficiently small such that if
\begin{itemize}
\item[(i)]   $ | a_0 |_{\H^1(\real)} + |a_0 |_{\L^1(\real)} < \epsilon_0$
and  
$ | b_0 |_{\H^2(\real)} + |b_0 |_{\L^1(\real)} < \epsilon_0$
\item[(ii)] $ | x^2 a_0 |_{\L^2(\real)} + | x^2 b_0 |_{\L^2(\real)} < \infty$,
\end{itemize}
then \reff{eqn:p-system} has a unique (mild) solution with
initial conditions $a_0$ and $b_0$. Moreover, there exist
functions $\{ g_{n}^{\pm} \}_{n=0}^{N}$ (independent of initial
conditions) and constants $C_N$, $\{d_{n}^{\pm}\}_{n=1}^{N}$ determined by the
initial conditions such that if we define
\begin{equs}
u(x,t) &= a(x-t,t) + b(x-t,t) ~~~\mbox{ and }~~~ v(x,t) = a(x+t,t) -  b(x+t,t)
\end{equs}
then
\begin{equa}[2]\label{eqn:asymptoticexpansion}
u(x,t) &= \frac{1}{\sqrt{1+t}}
g_0^{+}({\textstyle\frac{x}{\sqrt{1+t}}})
+ \sum_{n=1}^N \frac{1}{(1+t)^{1-\frac{1}{2^{n+1}}}} 
d_n^{+} g_{n}^{+}({\textstyle\frac{x}{\sqrt{1+t}}}) + R_u^N(x,t)
\\
v(x,t) &= \frac{1}{\sqrt{1+t}}
g_0^{-}({\textstyle\frac{x}{\sqrt{1+t}}})
+ \sum_{n=1}^N \frac{1}{(1+t)^{1-\frac{1}{2^{n+1}}}} 
d_n^{-} g_{n}^{-} ({\textstyle\frac{x}{\sqrt{1+t}}}) 
 + R_v^N(x,t)~,
\end{equa}
where the remainders $R_u^N$ and $R_v^N$ satisfy the estimates
\begin{equa}[2]\label{eqn:remainderestimates}
\sup_{t\geq0}(1+t)^{\frac{3}{4}-\frac{1}{2^{N+2}}}
\|R_{\{u,v\}}^N(\cdot,t)\|_{\L^2(\real)}&\leq C_N\\
\sup_{t\geq0}
(1+t)^{\frac{5}{4}-\frac{1}{2^{N+2}}}
\|\partial_x R_{\{u,v\}}^N(\cdot,t)\|_{\L^2(\real)} 
&\le C_N~.
\end{equa}
Furthermore, for $n\geq1$, the functions $g_n^{\pm}$ satisfy
$g_n^{\pm}(z)\sim |z|^{-1+2^{-n-1}}$ as $z\to\pm\infty$.
\end{theorem}
There is a slight incongruity in this result in that the norm in
which we estimate the remainder term is weaker than that we use
on the initial data; namely, we do not give
estimates for the remainder in $\H^2(\real)$, or in the
localization norms $\L^1(\real)$ and the weighted $\L^2(\real)$-norm
(on that aspect of the problem, see Remark \ref{rem:onestbete} below). Theorem
\ref{thm:maintheorem} actually holds for slightly more general
initial conditions than those satisfying (i)-(ii). Furthermore,
we will prove that the estimates (\ref{eqn:remainderestimates})
hold for all initial conditions $(a_0,b_0)$ in a subset ${\cal
D}_2\subset\H_1\times\H_2$ that is {\em positively invariant} under
the flow of \reff{eqn:p-system}. However, since the topology used
to define the subset ${\cal D}_2$ is somewhat non-standard, we have
chosen to state the result initially in this slightly weaker, but
hopefully more comprehensible, form to keep the introduction as
simple as possible.

\begin{remark}
\label{rem:onestbete}
It is interesting to note (see Proposition
\ref{prop:weightednorm} below) that
$\|x^2a(\cdot,t)\|_{\L^2(\real)}+\|x^2b(\cdot,t)\|_{\L^2(\real)}$
is finite for all finite $t>0$, but that the terms with $n\geq1$
in the asymptotic expansion do not satisfy this property due to
the long tails of the functions $g^{\pm}_{n}$. 
\end{remark}

\begin{remark}
As the asymmetry in the degree of $x$ derivatives in
(\ref{eqn:p-system}) suggests, we require more spatial regularity
from the second component (the $b$ variable) than from the first
(the $a$ variable). It is then natural to expect that $R_u^N$ or
$R_v^N$ are not necessarily in $\H^2$, but that only their
difference is.
\end{remark}

We conclude this section with a few remarks. Define
$u_{\pm}(x,t)=a(x,t)\pm b(x,t)$. Then the asymptotics of the
solutions of (\ref{eqn:p-system}) in the variables $u_{\pm}$ are
the same as those of the two dimensional (generalized) Burger's
equation
\begin{equa}[2]\label{eqn:modelproblem}
\partial_t u_{+}&=\partial_x^2u_{+}+\partial_xu_{+}
+\partial_x(c_{+}u_{+}^2-c_{-}u_{-}^2)\\
\partial_t u_{-}&=\partial_x^2u_{-}-\partial_xu_{-}
+\partial_x(c_{-}u_{-}^2-c_{+}u_{+}^2)~,
\end{equa}
where the constants $c_{\pm}$ are determined by the Hessian of
$g(a,b)$ at $a=b=0$ through
\begin{equs}
c_{\pm}=
\pm
\frac{1}{8}
(1,\pm1)\cdot
\myl{14}
\begin{matrix}
\partial_a^2g  & \partial_a\partial_bg\\
\partial_a\partial_bg & \partial_b^2g
\end{matrix}
\myr{14}
\my{|}{14}_{a=b=0}
\cdot
\vector{1}{\pm1}~.
\end{equs}
We will see that the hyperbolic effects manifest themselves
through the `source' terms $-c_{-}u_{-}^2$, respectively
$c_{+}u_{+}^2$ in the first, respectively second equation in
(\ref{eqn:modelproblem}). In particular, none of the terms
$g_n^{\pm}$ with $n\geq1$ would be present in the asymptotic
expansion if those terms were absent.

Finally, note that we have chosen to state Theorem
\ref{thm:maintheorem} for finite $N$. As it turns out, the sums
appearing in (\ref{eqn:asymptoticexpansion}) converge in the
limit as $N\to\infty$, in which case the estimates
(\ref{eqn:remainderestimates}) hold with time weights replaced by
$(1+t)^{\frac{3}{4}}\ln(2+t)^{-1}$ and
$(1+t)^{\frac{5}{4}}\ln(2+t)^{-1}$. The proof can easily be done
with the techniques used in this paper and is left to the reader.

The remainder of the paper is organized as follows: In Section
\ref{sec:cauchy}, we discuss the well-posedness of the Cauchy
problem (\ref{eqn:p-system}) in an appropriately defined
topology. In Section \ref{sec:asym}, we explain our strategy for
proving our main result, Theorem \ref{thm:maintheorem}, on the
long time asymptotics of solutions of (\ref{eqn:p-system}).
Namely, we decompose that proof into a series of simpler
sub-problems which are then tackled in subsequent sections: in
Sections \ref{sect:burgers} and \ref{sect:inhomogeneousheat}, we
investigate properties of solutions of Burger's type equations,
respectively of inhomogeneous heat equations, as they occur
naturally in the asymptotic analysis. In Section
\ref{sect:cauchyproof}, we collect some estimates that are used in
the proof of the well-posedness of (\ref{eqn:p-system}). Finally,
in Section \ref{sect:remainderestimates}, we specify the sense in
which the semigroup of the linearization of (\ref{eqn:p-system})
is close to heat kernels translating along the characteristics,
and we give estimates on the remainder terms occurring in Theorem
\ref{thm:maintheorem}.

\section{Cauchy problem}
\label{sec:cauchy}

To motivate our technical treatment of the problem and in particular
our choice of function spaces, we first note that upon taking
the Fourier transform of the linearization of
(\ref{eqn:p-system}), it follows that
\begin{equs}
\partial_t\vector{a}{b}=
\L\vector{a}{b}\equiv
\myl{14}
\begin{matrix}
0  & ik\\
ik &-2k^2
\end{matrix}
\myr{14}
\vector{a}{b}~.
\label{eqn:linearfourier}
\end{equs}
We then find that the (Fourier transform of) the 
semigroup associated with (\ref{eqn:linearfourier}) is
\begin{equs}
\ed^{\L t}=
\ed^{-k^2t}
\myl{15}
\begin{matrix}
\cos(kt\Delta)+\frac{k}{\Delta}\sin(kt\Delta) & 
\frac{i}{\Delta}\sin(kt\Delta)\\
\frac{i}{\Delta}\sin(kt\Delta) & 
\cos(kt\Delta)-\frac{k}{\Delta}\sin(kt\Delta)
\end{matrix}
\myr{15}~,
\label{eqn:defeLt}
\end{equs}
where $\Delta=\sqrt{1-k^2}$. The most important fact about the
semigroup $\ed^{\L t}$ is that it is close to $\ed^{\L_0 t}$,
the semigroup associated with the problem
\begin{equs}
\partial_t\vector{u}{v} =\L_0\vector{u}{v} \equiv \myl{14}
\begin{matrix} \partial_x^2+\partial_x & 0\\ 0 &
\partial_x^2-\partial_x \end{matrix} \myr{14} \vector{u}{v} ~.
\label{eqn:linearfourieruv}
\end{equs}
Formally, $\ed^{\L_0 t}$ can be obtained by setting $\Delta=1$ in
$\ed^{\L t}$ and by conjugating with the matrix
\begin{equs}
\S\equiv
\myl{14}
\begin{matrix}
1 & 1\\
1 & -1
\end{matrix}
\myr{14}~.
\label{eqn:defS}
\end{equs}
These two operations correspond to a long wavelength expansion
and a change of dependent variables to quantities that move along
the characteristics. More precisely, we will prove that $\ed^{\L
t}$ satisfies the intertwining property 
\begin{equs}
\S\ed^{\L t}&\approx\ed^{\L_0 t}\S~,
\end{equs}
where the symbol $\approx$ means that the action of these two
operators is the same in the large scale -- long time limit; see
Lemma \ref{lem:closetoheat} at the beginning of Section
\ref{sect:remainderestimates} for details.

Furthermore, $\ed^{\L t}$ satisfies parabolic-like estimates
\begin{equs}
|\ed^{\L t}|&\leq
C\ed^{-\min(k^2,1)\frac{t}{4}}
\myl{15}
\begin{matrix}
1 & 
\frac{1}{\sqrt{1+k^2}}\\
\frac{1}{\sqrt{1+k^2}} & 
1
\end{matrix}
\myr{15}~,
\label{eqn:estimateonkernel}
\\
\my{|}{15}
\ed^{\L t}
\vector{0}{ik}
\my{|}{15}
&
\leq
C\frac{\ed^{-\min(k^2,1)\frac{t}{4}}}{\sqrt{t}}
\vector{1}{\frac{1}{\sqrt{1+k^2}}}
\label{eqn:estimateonDkernel}
\end{equs}
uniformly in $t\geq0$ and $k\in{\bf R}$.

Hence, to summarize, $\ed^{\L t}$ behaves like a superposition of
heat kernels translating along the characteristics of the
underlying hyperbolic problem. In view of the above observations
as well as of classical techniques for parabolic PDE's, see e.g.
\cite{temam:1997,bricmont:1994}, we will consider
(\ref{eqn:p-system}) in the following (somewhat non-standard)
topology (cf also \cite{gvb:2006}):
\begin{definition}
\label{def:functionspaces}
We define ${\cal B}_0$, resp. ${\cal B}$, as the closure of 
${\cal C}_0^{\infty}({\bf R},{\bf R}^2)$, resp. ${\cal C}_0^{\infty}({\bf
R}\times[0,\infty),{\bf R}^2)$, under the norm $|\cdot|$, resp.
$\|\cdot\|$, where for ${\bf z}_0=(a_0,b_0):{\bf R}\to{\bf R}^2$
and ${\bf z}=(a,b):{\bf R}\times[0,\infty)\to{\bf R}^2$, we define
\begin{equs}
|{\bf z}_0|=
 \|\widehat{\bf z}_0\|_{\infty}
+\|{\bf z}_0\|_2
+\|\D{\bf z}_0\|_2
+\|\D^2b_0\|_2~,~~~
\|{\bf z}\|&=
\|\hat{\bf z}\|_{\infty,0}
+\|{\bf z}\|_{2,\frac{1}{4}}
+\|\D{\bf z}\|_{2,\frac{3}{4}}+
\|\D^2b\|_{2,\frac{5}{4}^{\star}}~.
\end{equs}
Here $(Da)(x,t)\equiv\partial_xa(x,t)$, $\hat{a}(k,t)$ is the
Fourier transform of $a(x,t)$,
\begin{equs}
\|f\|_{p,q}=\sup_{t\geq0}(1+t)^q\|f(\cdot,t)\|_p~,~~~~
\|f\|_{p,q^{\star}}=\sup_{t\geq0}\frac{(1+t)^q}{\ln(2+t)}\|f(\cdot,t)\|_p
\end{equs}
and $\|\cdot\|_p$ is the standard $\L^p({\bf R})$ norm.
\end{definition}
Before turning to the Cauchy problem with initial data in ${\cal
B}_0$ we collect a few comments on our choice of function spaces.

Consider first the requirements on the initial conditions in
(\ref{eqn:p-system}). While the use of $\H^1$ space is quite
natural in this context, we choose to replace the $\L^1$ norm by
the (weaker) control of the $\L^{\infty}$ norm in Fourier space.
This has the great advantage that all estimates can then be done
in Fourier space, where the semigroup $\ed^{\L t}$ has the
simple, explicit, form (\ref{eqn:defeLt}).

In turn, our choice of $q$-exponents in the norm $\|\cdot\|$ is
motivated by the fact that these are the highest possible
exponents for which the $\|\cdot\|$-norm of the leading order
asymptotic term $\frac{1}{\sqrt{1+t}}g_0(\frac{x}{\sqrt{1+t}})$
is bounded. Note also that for the linear evolution
(\ref{eqn:linearfourier}), we have
\begin{equs}
\|\ed^{\L t}{\bf z}_0\|\leq
C|{\bf z}_0|~,
\label{eqn:linearestimate}
\end{equs}
since $\hat{j}(k,t)=\ed^{-\min(k^2,1)t}u_0(k)$ satisfies
\begin{equs}
\|\D^nj(\cdot,t)\|_{2}
\leq C
\myl{12}
\ed^{-t}\|\D^nu_0\|_2
+\min\myl{10}
t^{-\frac{1}{4}-\frac{n}{2}}\|\hat{u}_0\|_{\infty},
\|D^nu_0\|_2
\myr{10}
\myr{12}
\end{equs}
for all $n=0,1,\ldots$.

Finally, we note that for admissible nonlinearities in the sense
of Definition \ref{def:admissible}, the map
$h(a,b)=f(a,b)\partial_xb+g(a,b)=h({\bf z})$ satisfies
\begin{equs}
\label{eqn:H}
\|h({\bf z})\|_{1,\frac{1}{2}}+
\|h({\bf z})\|_{2,\frac{3}{4}}+
\|Dh({\bf z})\|_{2,\frac{5}{4}}&\leq 
C\|{\bf z}\|^2~,\\
\|h({\bf z}_1)-h({\bf z}_2)\|_{1,\frac{1}{2}}+
\|h({\bf z}_1)-h({\bf z}_2)\|_{2,\frac{3}{4}}&\leq 
C\|{\bf z}_1-{\bf z}_2\|(\|{\bf z}_1\|+\|{\bf z}_2\|)~,
\label{eqn:DHone}
\\
\|D(h({\bf z}_1)-h({\bf z}_2))\|_{2,\frac{5}{4}}&\leq 
C\|{\bf z}_1-{\bf z}_2\|(\|{\bf z}_1\|+\|{\bf z}_2\|)~.
\label{eqn:DHtwo}
\end{equs}
We are now fully equipped to study the Cauchy problem (\ref{eqn:p-system})
in ${\cal B}$:
\begin{theorem}\label{thm:cauchy}
For all ${\bf z}_0\in{\cal B}_0$ with $|{\bf
z}_0|=|(a_0,b_0)|\leq\epsilon_0$ small enough,
the Cauchy problem (\ref{eqn:p-system})
is (locally) well posed in ${\cal B}$ if the nonlinearities are
admissible in the sense of Definition \ref{def:admissible}.
In particular, the solution satisfies $\|{\bf z}\|\leq
c\epsilon_0$ for some $c>1$ and is unique among functions in
${\cal B}$ satisfying this bound.
\end{theorem}

\begin{proof}
Upon taking the Fourier transform of (\ref{eqn:p-system}), we get
\begin{equs}
\partial_t\vector{a}{b}=
\myl{14}
\begin{matrix}
0  & ik\\
ik &-2k^2
\end{matrix}
\myr{14}
\vector{a}{b}+
\vector{0}{ikh}~,
\label{eqn:fourier}
\end{equs}
which gives the following representation for the solution
\begin{equs}
{\bf z}(t)\equiv
\vector{a(t)}{b(t)}
=\ed^{\L t}\vector{a_0}{b_0}
+\int_0^t
\hspace{-2mm}
{\rm d}s~\ed^{\L(t-s)}
\vector{0}{\partial_xh({\bf z}(s))}
\equiv
\ed^{\L t}{\bf z}_0+
{\cal N}[{\bf z}](t)
~.
\label{eqn:defN}
\end{equs}
We will prove below that for all ${\bf z}_i\in{\cal B}$, $i=1,2$,
we have
\begin{equs}
\|{\cal N}[{\bf z}]\|\leq C\|{\bf z}\|^2~~\mbox{ and }~~
\|{\cal N}[{\bf z_1}]-{\cal N}[{\bf z_2}]\|\leq 
C\|{\bf z}_1-{\bf z}_2\|(\|{\bf z}_1\|+\|{\bf z}_2\|)
\label{eqn:contract}
\end{equs}
for some constant $C$. The proof of Theorem \ref{thm:cauchy}
then follows from the fact that for all ${\bf z}_0\in{\cal B}_0$
with $|{\bf z}_0|\leq \epsilon_0$ small enough and $c>1$, the
r.h.s. of (\ref{eqn:defN}) defines a contraction map from some
(small) ball of radius $c\epsilon_0$ in ${\cal B}$ onto itself.

The general rule for proving the various estimates involved in
(\ref{eqn:contract}) is to split the integration interval into
two parts, with $s\in{\cal I}_1\equiv[0,\frac{t}{2}]$ and
$s\in{\cal I}_2\equiv[\frac{t}{2},t]$. In ${\cal I}_1$, we place
as many derivatives (or equivalently, factors of $k$) as possible
on the semigroup $\ed^{\L(t-s)}$, while on ${\cal I}_2$, (most
of) these derivatives need to act on $h$, since the
integral would otherwise be divergent at $s=t$.

Additional difficulties arise from the fact that $\ed^{\L t}$ has
very little smoothing properties (slow or no decay in $k$ as
$|k|\to\infty$), so that in some cases we need to consider
separately the large-$k$ part and the small-$k$ part of the
$\L^2$ norm, say. This is done through the use of $\P$, defined
as the Fourier multiplier with the characteristic function on
$[-1,1]$.

We decompose the proof of $\|{\cal N}[{\bf z}]\|\leq C\|{\bf
z}\|^2$ into that of
\begin{equs}
\|{\cal N}[{\bf z}]\|&\leq
 \|\widehat{{\cal N}[{\bf z}]}\|_{\infty,0}
+\|{\cal N}[{\bf z}]\|_{2,\frac{1}{4}}
+\|\P \D{\cal N}[{\bf z}]\|_{2,\frac{3}{4}}
+\|(1-\P)\D{\cal N}[{\bf z}]\|_{2,\frac{3}{4}}\\
&\phantom{=~}
+\|(1-\P)\D^2{\cal N}[{\bf z}]_2\|_{2,\frac{5}{4}^{\star}}
+\|(1-\Q)\P\D^2{\cal N}[{\bf z}]_2\|_{2,\frac{5}{4}^{\star}}
+\|\Q\P\D^2{\cal N}[{\bf z}]_2\|_{2,\frac{5}{4}^{\star}}
\\&\leq C
\|{\bf z}\|^2
~,
\label{eqn:split}
\end{equs}
where $\Q$ is the characteristic function for $t\geq1$ and ${\cal
N}[{\bf z}]_2$ denotes the second component of ${\cal N}[{\bf
z}]$.

We now consider $\|\P\D{\cal N}[{\bf z}]\|_{2,\frac{3}{4}}$ as an
example of the way we prove the above estimates. We have
\begin{equs}
\|\P\D{\cal N}[{\bf z}](\cdot,t)\|_2
&\leq
\|h({\bf z})\|_{2,\frac{3}{4}}
\my{(}{13}\sup_{|k|\leq1,\tau\geq0}
\hspace{-3mm}
|k|\sqrt{\tau}\ed^{-\frac{k^2\tau}{4}}
\my{)}{13}
\int_0^{\frac{t}{2}}
\hspace{-3mm}{\rm d}s~
\frac{(1+s)^{-\frac{3}{4}}}{t-s}
\\&\phantom{=~}
+
\|\D h({\bf z})\|_{2,\frac{5}{4}}
\my{(}{13}
\sup_{|k|\leq1,\tau\geq0}
\hspace{-3mm}
\ed^{-\frac{k^2\tau}{4}}
\my{)}{13}
\int_{\frac{t}{2}}^t
\hspace{-2mm}{\rm d}s~
\frac{(1+s)^{-\frac{5}{4}}}{\sqrt{t-s}}
\\&\leq
C\|{\bf z}\|^2
\myl{12}
\frac{2}{t}
\int_0^{\frac{t}{2}}
\hspace{-3mm}
\frac{{\rm d}s~}{(1+s)^{\frac{3}{4}}}+
\frac{1}{(1+\frac{t}{2})^{\frac{5}{4}}}
\int_{\frac{t}{2}}^t
\hspace{-2mm}
\frac{{\rm d}s~}{\sqrt{t-s}}
\myr{12}
\leq C\|{\bf z}\|^2 (1+t)^{-\frac{3}{4}}
\label{eqn:firstB}
\end{equs}
for all $t\geq0$, which shows that $\|\P\D{\cal N}[{\bf
z}]\|_{2,\frac{3}{4}}\leq C\|{\bf z}\|^2$. All other estimates in
(\ref{eqn:split}) can be done similarly; we postpone their proof
to Section \ref{sect:cauchyproof} below.

Finally, we note that the Lipschitz-type estimate in
(\ref{eqn:contract}) can be obtained in the same manner, {\em
mutatis mutandis}, due to the similarity between
(\ref{eqn:DHone}) and (\ref{eqn:DHtwo}) with (\ref{eqn:H}); we
omit the details.
\end{proof} 

We can now turn to the question of the asymptotic structure of
the solutions of (\ref{eqn:p-system}) provided by Theorem
\ref{thm:cauchy}.
Note that already if we wanted to prove that $\ed^{\L t}{\bf
z}_0$ satisfies `Gaussian asymptotics' we would need more
localization properties on ${\bf z}_0$ than those provided by the
${\cal B}_0$-topology. It will turn out to be sufficient to
require ${\bf z}_0\in{\cal B}_0\cap \L^2(\real,x^m\d x)$ for
(some) $m\geq2$. We now prove that this requirement is
{\em forward invariant} under the flow of (\ref{eqn:p-system}):
\begin{proposition}\label{prop:weightednorm}
Let $\rho_m(x)=|x|^m$ and define
\begin{equs}
{\cal D}_m=\my{\{}{12}{\bf z}_0\in{\cal B}_0\mbox{ such that }|{\bf
z}_0|+\|\rho_m{\bf z}_0\|_2<\infty
\my{\}}{12}~.
\end{equs}
If ${\bf z}_0\in{\cal D}_m$ and $|{\bf z}_0|\leq\epsilon_0$ such
that Theorem \ref{thm:cauchy} holds, then the corresponding
solution ${\bf z}(t)$ of (\ref{eqn:p-system}) satisfies ${\bf
z}(t)\in{\cal D}_m$ for all finite $t>0$. Furthermore, there holds 
$|{\bf z}(t)|\leq (1+\delta)\epsilon_0$ for some (small) constant
$\delta$.
\end{proposition}

\begin{proof}
Note first that by Theorem \ref{thm:cauchy}, 
$|{\bf z}(t)|\leq\|{\bf z}\|\leq(1+\delta)\epsilon_0$ since 
${\bf z}_0\in{\cal B}_0$ and $|{\bf z}_0|\leq\epsilon_0$. Then,
fix $m\in{\bf N}$, $m\geq1$. The proof of Theorem
\ref{thm:cauchy} can easily be adapted to show that
(\ref{eqn:p-system}) is {\em locally} (in time) well posed in
${\cal D}_m$. Global existence then follows from the fact that the quantity
\begin{equs}
N(t)=\frac{1}{2}
\|\rho_m{\bf z}(\cdot,t)\|^2
=\frac{1}{2}\int_{-\infty}^{\infty}
\hspace{-4mm}{\rm d}x~
|x|^m(a(x,t)^2+b(x,t)^2)
\end{equs}
grows {\em at most exponentially} as $t\to\infty$.
Namely, we have
\begin{equs}
\partial_t N(t)&=
\int_{-\infty}^{\infty}
\hspace{-4mm}{\rm d}x~
|x|^m\myl{12}
\partial_x(ab)+
2b\partial_x^2b
+b\partial_x\myl{10}f(a,b)\partial_xb+g(a,b)\myr{10}
\myr{12}\\
&=-
\int_{-\infty}^{\infty}
\hspace{-4mm}{\rm d}x~
m|x|^{m-1}\sign(x)
\myl{12}
b(a+g(a,b))+(2+f(a,b))b\partial_xb
\myr{12}
\\&\phantom{=~}
-\int_{-\infty}^{\infty}
\hspace{-4mm}{\rm d}x~
|x|^m
(\partial_xb)^2
\myl{10}
2+f(a,b)
\myr{10}
\\&\leq
\int_{-\infty}^{\infty}
\hspace{-4mm}{\rm d}x~
\myl{10}
(m-1)^{m-1}+|x|^m
\myr{10}
\my{|}{12}
b(a+g(a,b))+(2+f(a,b))b\partial_xb
\my{|}{12}
\\&\phantom{=~}
-\int_{-\infty}^{\infty}
\hspace{-4mm}{\rm d}x~
|x|^m
(\partial_xb)^2
\myl{10}
2+f(a,b)
\myr{10}
\\&\leq
\int_{-\infty}^{\infty}
\hspace{-4mm}{\rm d}x~
\myl{10}
(m-1)^{m-1}+|x|^m
\myr{10}
\myl{12}
|b(a+g(a,b))|+
2^{-1}|2+f(a,b)|b^2
\myr{12}
\\&\leq
C_1(m,\epsilon_0)+C_2(\epsilon_0)N(t)~,
\end{equs}
due to the estimates $\|f(a,b)\|_{\infty}\leq C\epsilon_0\ll2$ and
$\|\frac{g(a,b)}{\sqrt{a^2+b^2}}\|_{\infty}\leq C\epsilon_0$.
\end{proof}

\section{Asymptotic structure - Proof of Theorem \ref{thm:maintheorem}}
\label{sec:asym}

We can now state our main result on the asymptotic structure of
solutions of (\ref{eqn:p-system}) in a definitive manner:
\begin{theorem}\label{thm:asymptoticsrestated}
Let ${\cal D}_m$ be as in Proposition \ref{prop:weightednorm}
with $m\geq2$, let ${\bf z}_0\in{\cal D}_m$ with $|{\bf
z}_0|\leq\epsilon_0$ such that Theorem \ref{thm:cauchy} holds and
define
\begin{equs}
u(x,t) &= a(x-t,t) + b(x-t,t) ~~~\mbox{ and }~~~
v(x,t) = a(x+t,t) -  b(x+t,t)
\end{equs}
for the corresponding solution ${\bf z}(t)=(a(t),b(t))$ of
(\ref{eqn:p-system}). Then there exist functions $\{ g_{n}^{\pm}
\}_{n=0}^{N}$ (independent of ${\bf z}_0$) and constants $C_N$,
$\{d_{n}^{\pm}\}_{n=1}^N$ determined by ${\bf z}_0$ such that
\begin{equa}[2]\label{eqn:asymptoticexpansionrap}
u(x,t) &= \frac{1}{\sqrt{1+t}}
g_0^{+}({\textstyle\frac{x}{\sqrt{1+t}}})
+ \sum_{n=1}^N \frac{1}{(1+t)^{1-\frac{1}{2^{n+1}}}} 
d_n^{+} g_{n}^{+}({\textstyle\frac{x}{\sqrt{1+t}}}) + R_u^N(x,t)
\\
v(x,t) &= \frac{1}{\sqrt{1+t}}
g_0^{-}({\textstyle\frac{x}{\sqrt{1+t}}})
+ \sum_{n=1}^N \frac{1}{(1+t)^{1-\frac{1}{2^{n+1}}}} 
d_n^{-} g_{n}^{-} ({\textstyle\frac{x}{\sqrt{1+t}}}) 
 + R_v^N(x,t)~,
\end{equa}
where the remainders $R_u^N$ and $R_v^N$ satisfy the estimates
\begin{equa}[2]\label{eqn:remainderestimatesrap}
\sup_{t\geq0}(1+t)^{\frac{3}{4}-\frac{1}{2^{N+2}}}
\|R_{\{u,v\}}^N(\cdot,t)\|_{\L^2(\real)}&\leq C_N\\
\sup_{t\geq0}
(1+t)^{\frac{5}{4}-\frac{1}{2^{N+2}}}
\|\partial_x R_{\{u,v\}}^N(\cdot,t)\|_{\L^2(\real)} 
&\le C_N~.
\end{equa}
Furthermore, for $n\geq1$, the functions $g_n^{\pm}$ satisfy
$g_n^{\pm}(z)\sim |z|^{-1+2^{-n-1}}$ as $z\to\pm\infty$.
\end{theorem}

\begin{remark}
As will be apparent from the proof of Theorem
\ref{thm:asymptoticsrestated}, any hyperbolic-parabolic
system of the form
\begin{equs}
\partial_t{\bf z}+f({\bf z})_x=(B({\bf z}){\bf z}_x)_x
\end{equs}
with admissible nonlinearities in the sense of (the natural
extension of) Definition \ref{def:admissible} gives rise to
solutions having the same asymptotic structure as those of the
p-system as long as the following two conditions are satisfied:
\begin{enumerate}
\item\label{item:inter} There exist two matrices ${\cal S}$ and
${\rm A}$ with ${\cal S}$ non-singular and ${\rm A}$ diagonal
having eigenvalues of multiplicity $1$ for which ${\cal S}\ed^{\L
t}\approx\ed^{\L_0 t}{\cal S}$ in the sense of Lemma
\ref{lem:closetoheat} (see Section
\ref{sect:remainderestimates}), where
$\L_0=\partial_x^2+{\rm A}\partial_x$ and
$\L=B(0)\partial_x^2-f'(0)\partial_x$.
\item\label{item:cauchy} The Cauchy problem with initial
condition in the corresponding functional space (the natural
extension of ${\cal B}_0$ to the problem considered) is well
posed and satisfies the analogues of Theorem \ref{thm:cauchy} and
Proposition \ref{prop:weightednorm}.
\end{enumerate}
\end{remark}

We now briefly comment on the above assumptions for specific
systems such as the `full gas dynamics' and the MHD system. The
intertwining property of item \ref{item:inter} above is proved in
\cite{liu:1997} for quite general systems, though not in exactly
the same topology as that used in Lemma \ref{lem:closetoheat}. As
for item \ref{item:cauchy}, local well-posedness for initial data
in ${\cal B}_0$ is certainly not an issue, the only difficulty is
to prove that the various norms of Definition
\ref{def:functionspaces} exhibit `parabolic-like' decay as
$t\to\infty$. This is very likely to hold, particularly for
systems satisfying item \ref{item:inter}.

While the variables $(a,b)$ are adapted to the study of the
Cauchy problem because of the inherent asymmetry of spatial
regularity in (\ref{eqn:p-system}), they are not the best
framework for studying the asymptotic structure of the solutions
to (\ref{eqn:p-system}). It turns out to be more convenient to
change variables to quantities that move along the
characteristics. We thus define
\begin{equs}
\vector{u(x,t)}{v(x,t)}
\equiv
\myl{14}
\begin{matrix}
\T^{-1} & 0\\
0 & \T
\end{matrix}
\myr{14}
\myl{14}
\begin{matrix}
1 & 1\\
1 & -1
\end{matrix}
\myr{14}
\vector{a(x,t)}{b(x,t)}
\equiv 
\myl{14}
\begin{matrix}
\T^{-1} & 0\\
0 & \T
\end{matrix}
\myr{14}\S{\bf z}(x,t)
~,
\end{equs}
where $\T$ is the translation operator defined by
\begin{equs}
(\T f)(x,t)=f(x+t,t)~~\mbox{ or equivalently by }~~
\widehat{\T f}(k,t)=\ed^{ikt}\hat{f}(k,t)~.
\label{eqn:translationdef}
\end{equs}
Note in passing that
\begin{equs}
a(x,t)=\frac{1}{2}\myl{14}
u(x+t,t)+v(x-t,t)
\myr{14}~~\mbox{ and }~~
b(x,t)=\frac{1}{2}\myl{14}
u(x+t,t)-v(x-t,t)
\myr{14}~.
\end{equs}
We then use the fact that ${\bf z}$ satisfies the integral
equation
\begin{equs}
\S{\bf z}(t)&=\S\ed^{\L t}{\bf z}_0+
\int_0^t
\hspace{-2mm}
{\rm d}s~\S\ed^{\L(t-s)}
\vector{0}{\partial_xh({\bf z}(s))}\\
&=\ed^{\L_0 t}\S{\bf z}_0+
\int_0^t
\hspace{-2mm}
{\rm d}s~
\ed^{\L_0(t-s)}\S~
\vector{0}{\partial_xg_0({\bf z}(s))}
+{\cal R}[{\bf z}](t)~,
\label{eqn:representation}
\end{equs}
where
\begin{equs}
{\cal R}[{\bf z}](t)&=
\myl{10}
\S\ed^{\L t}-\ed^{\L_0 t}\S
\myr{10}{\bf z}_0
+
\int_0^t
\hspace{-2mm}
{\rm d}s~
\my{[}{14}
\S\ed^{\L(t-s)}
\vector{0}{\partial_xh({\bf z}(s))}
-
\ed^{\L_0(t-s)}\S~
\vector{0}{\partial_xg_0({\bf z}(s))}
\my{]}{14}~.
\end{equs}
To justify the notation, which suggests that ${\cal R}$ is a
remainder term, we will prove in Section
\ref{sect:remainderestimates} that ${\cal R}[{\bf z}]
=({\cal R}_u[{\bf z}],{\cal R}_v[{\bf z}])$ satisfies
the improved decay rates
\begin{equs}
\|{\cal R}_{\{u,v\}}[{\bf z}]\|_{2,\frac{3}{4}^{\star}}+
\|\D{\cal R}_{\{u,v\}}[{\bf z}]\|_{2,\frac{5}{4}^{\star}}\leq C\epsilon_0~,
\label{eqn:onRannounce}
\end{equs}
because of the intertwining relation $\S\ed^{\L
t}\approx\ed^{\L_0 t}\S$ (see Lemma \ref{lem:closetoheat}) and
the fact that $h({\bf z})=g_0({\bf z})+h.o.t.$.

Recalling that $g_0$ is quadratic (cf Definition
\ref{def:admissible}), we will write
\begin{equs}
g_0({\bf z})&=c_{+}(a+b)^2-c_{-}(a-b)^2+c_{3}(a+b)(a-b)\\
&=c_{+}(\T u)^2-c_{-}(\T^{-1}v)^2
+c_{3}(\T u)(\T^{-1}v)
\end{equs}
for ${\bf z}=(a,b)$. We thus find from (\ref{eqn:representation})
that $u$ and $v$ satisfy
\begin{equs}
u(t)&=\ed^{\partial_x^2t}(a_0+b_0)
+\partial_x\int_0^{t}
\hspace{-2mm}{\rm d}s~\ed^{\partial_x^2(t-s)}
\myl{12}
 c_{+}u(s)^2
-c_{-}\T^{-2}v(s)^2
\myr{12}
\\
&\phantom{=~}
+\T^{-1}{\cal R}_{u}[{\bf z}](t)
+c_{3}\partial_x\int_0^{t}
\hspace{-2mm}{\rm d}s~\ed^{\partial_x^2(t-s)}
\T^{-1}\myl{10}
(\T u(s))(\T^{-1}v(s))
\myr{10}
\label{eqn:firstintegu}
~,\\
v(t)&=\ed^{\partial_x^2t}(a_0 - b_0)
+\partial_x\int_0^{t}
\hspace{-2mm}{\rm d}s~\ed^{\partial_x^2(t-s)}
\myl{12}
 c_{-}v(s)^2
-c_{+}\T^{2}u(s)^2
\myr{12}
\\
&\phantom{=~}
+\T{\cal R}_{v}[{\bf z}](t)
-c_{3}\partial_x\int_0^{t}
\hspace{-2mm}{\rm d}s~\ed^{\partial_x^2(t-s)}
\T\myl{10}
(\T u(s))(\T^{-1}v(s))
\myr{10}
\label{eqn:firstintegv}
~.
\end{equs}
Note that, but for the presence of the second lines in
(\ref{eqn:firstintegu}) and (\ref{eqn:firstintegv}), these
expressions are precisely Duhamel's formula for the solution of
the model problem (\ref{eqn:modelproblem}), written in terms of
$u=\T^{-1} u_{+}$ and $v=\T u_{-}$. The next step is to write
\begin{equs}
u=u_{\star}+R_u^N=u_0+u_1+R_u^N~~\mbox{
and }~~
v=v_{\star}+R_v^N=v_0+v_1+R_v^N~,
\end{equs}
considering $R_u^N$ and $R_v^N$ as new `unknowns' and
\begin{equa}[2]
u_0(x,t)&=\frac{1}{\sqrt{1+t}}
g_0^{+}({\textstyle\frac{x}{\sqrt{1+t}}})~,&
u_1(x,t)&=
\sum_{n=1}^{N}
\frac{1}{(1+t)^{1-\frac{1}{2^{n+1}}}}
d_n^{+}g^{+}_{n}({\textstyle\frac{x}{\sqrt{1+t}}})
\\
v_0(x,t)&=\frac{1}{\sqrt{1+t}}
g_0^{-}({\textstyle\frac{x}{\sqrt{1+t}}})~&~\mbox{ and }~~
v_1(x,t)&=
\sum_{n=1}^{N}
\frac{1}{(1+t)^{1-\frac{1}{2^{n+1}}}}
d_n^{-}g^{-}_{n}({\textstyle\frac{x}{\sqrt{1+t}}})
\label{eqn:defu0u1}
\end{equa}
for some coefficients $\{d_n^{\pm}\}_{n=1}^{N}$ and functions
$\{g_{n}^{\pm}\}_{n=0}^{N}$ to be determined later.

We now use
\begin{equs}
u^2&=(u-u_{\star})(u+u_{\star})+u_{\star}^2=
R_u^N(u+u_{\star})+u_1^2+2u_0u_1+u_0^2~,\\
v^2&=(v-v_{\star})(v+v_{\star})+v_{\star}^2=
R_v^N(v+v_{\star})+v_1^2+2v_0v_1+v_0^2~,\\
(\T u)(\T^{-1}v)&=
 (\T R_u^N)\T^{-1}\myl{12}\frac{v+v_{\star}}{2}\myr{12}
+(\T^{-1}R_v^N)\T\myl{12}\frac{u+u_{\star}}{2}\myr{12}
+(\T u_{\star})(\T^{-1}v_{\star})~.
\end{equs}

Since
\begin{equs}[2]
g_0^{+}(x)&=u_0(x,0)~,&~~~
u_1(x,0)&=\sum_{n=1}^{N}d^{+}_ng_{n}^{+}(x)~,\\
g_0^{-}(x)&=v_0(x,0)~&~\mbox{ and }~~
v_1(x,0)&=\sum_{n=1}^{N}d^{-}_ng_{n}^{-}(x)~,
\end{equs}
we find that $R_u^N$ and $R_v^N$ satisfy
\begin{equs}
R_u^N(t)&=
\ed^{\partial_x^2t}(a_0+b_0-g_0^{+})
\\&\phantom{= }
+\my{[}{14}
\ed^{\partial_x^2t}u_0(0)
+ c_{+}\partial_x\int_0^{t}
\hspace{-2mm}{\rm d}s~\ed^{\partial_x^2(t-s)}
u_0(s)^2
\my{]}{14}-u_0(t)
\\&\phantom{= }
+\my{[}{14}
\ed^{\partial_x^2 t}u_1(0)
+ 2c_{+}\partial_x\int_0^{t}
\hspace{-2mm}{\rm d}s~\ed^{\partial_x^2(t-s)}
u_0(s)u_1(s)
\my{]}{14}-u_1(t)
\\&\phantom{= }
- c_{-}\my{[}{14}
\partial_x\int_0^{t}
\hspace{-2mm}{\rm d}s~\ed^{\partial_x^2(t-s)}
\T^{-2}
\myl{12}
(v_0(s)^2+2v_0(s)v_1(s))
\myr{12}
\my{]}{14}
-
\sum_{n=1}^{N}
\ed^{\partial_x^2 t}
d^{+}_ng_{n}^{+}
\\
&\phantom{= }
+\widetilde{\cal R}_{u}[{\bf z},{\bf R}^{N}](t)
+\T^{-1}{\cal R}_{u}[{\bf z}](t)~,
\label{eqn:tractable_u}\\[4mm]
R_v^N(t)&=
\ed^{\partial_x^2t}(a_0-b_0-g_0^{-})
\\&\phantom{= }
+\my{[}{14}
\ed^{\partial_x^2t}v_0(0)
+ c_{-}\partial_x\int_0^{t}
\hspace{-2mm}{\rm d}s~\ed^{\partial_x^2(t-s)}
v_0(s)^2
\my{]}{14}-v_0(t)
\\&\phantom{= }
+\my{[}{14}
\ed^{\partial_x^2 t}v_1(0)
+ 2c_{-}\partial_x\int_0^{t}
\hspace{-2mm}{\rm d}s~\ed^{\partial_x^2(t-s)}
v_0(s)v_1(s)
\my{]}{14}-v_1(t)
\\&\phantom{= }
-c_{+}\my{[}{14}
\partial_x\int_0^{t}
\hspace{-2mm}{\rm d}s~\ed^{\partial_x^2(t-s)}
\T^{2}
\myl{12}
(u_0(s)^2+2u_0(s)u_1(s))
\myr{12}
\my{]}{14}
-
\sum_{n=1}^{N}
\ed^{\partial_x^2 t}
d^{-}_ng_{n}^{-}
\\
&\phantom{= }
+\widetilde{\cal R}_{v}[{\bf z},{\bf R}^{N}](t)
+\T{\cal R}_{v}[{\bf z}](t)~,
\label{eqn:tractable_v}
\end{equs}
where
\begin{equs}
\widetilde{\cal R}_{u}[{\bf z},{\bf R}^{N}](t)&=
 c_{+}\E_0   [h_{1,u}+h_{3,u}](t)
-c_{-}\E_{-2}[h_{1,v}+h_{3,v}](t)
+c_{3}\E_{-1}[h_{2}  +h_{4}  ](t)
~,
\\
\widetilde{\cal R}_{v}[{\bf z},{\bf R}^{N}](t)&=
c_{-}\E_0   [h_{1,v}+h_{3,v}](t)
-c_{+}\E_{ 2}[h_{1,u}+h_{3,u}](t)
-c_{3}\E_{ 1}[h_{2}  +h_{4}  ](t)~,
\end{equs}
with ${\bf R}^N=(R_u^N,R_v^N)$, 
\begin{equs}[2]
\E_{\sigma}[h](t)&=
\partial_x\int_0^{t}
\hspace{-2mm}{\rm d}s~\ed^{\partial_x^2(t-s)}~
\T^{\sigma}h(s)~&~\mbox{and}\\
h_{1,u}&=R_u^N(u+u_{\star})~,~~~h_{3,u}=u_1^2~,&~~~
h_2&=(\T R_u^N)\T^{-1}\myl{12}\frac{v+v_{\star}}{2}\myr{12}
+(\T^{-1}R_v^N)\T\myl{12}\frac{u+u_{\star}}{2}\myr{12}\\
h_{1,v}&=R_v^N(v+v_{\star})~,~~~h_{3,v}=v_1^2~,&~~~
h_4&=(\T u_{\star})(\T^{-1}v_{\star})~.
\end{equs}
Note that we can write (\ref{eqn:tractable_u}) and
(\ref{eqn:tractable_v}) as ${\bf R}^N={\cal F}[{\bf z},{\bf
R}^N]$. If we now consider ${\bf z}$ fixed, we can interpret
${\bf R}^N={\cal F}[{\bf z},{\bf R}^N]$ as an equation for
${\bf R}^N$ which can be solved via a contraction mapping
argument. Namely, we will prove that if $\|{\bf
z}\|\leq C\epsilon_0$,  ${\bf R}^N\mapsto{\cal F}[{\bf z},{\bf
R}^N]$ defines a contraction map inside the ball
\begin{equs}
\|R_u^N\|_{2,\frac{3}{4}-\epsilon}+
\|\D R_u^N\|_{2,\frac{5}{4}-\epsilon}+
\|R_v^N\|_{2,\frac{3}{4}-\epsilon}+
\|\D R_v^N\|_{2,\frac{5}{4}-\epsilon}
\leq C
\label{eqn:restestim}
\end{equs}
for $\epsilon=2^{-N-2}$, provided $\{g^{\pm}_{n}\}_{n=0}^{N}$ and 
$\{d^{\pm}_n\}_{n=1}^{N}$ are appropriately chosen.

Basically, we will choose $u_0$, $v_0$, $u_1$ and $v_1$ in such a
way that the second and third lines of (\ref{eqn:tractable_u})
and (\ref{eqn:tractable_v}) vanish. Note that if, for instance,
we set the second, respectively third lines of
(\ref{eqn:tractable_u}) and (\ref{eqn:tractable_v}) equal to
zero, the resulting equalities are nothing but Duhamel's formulae
for Burger's equations for $u_0$ and $v_0$, respectively for
linearized Burger's equations for $u_1$ and $v_1$. Properties of
solutions to these types of equations are studied in detail in
Section \ref{sect:burgers} below.

Once $u_0$, $v_0$, $u_1$ and $v_1$ are fixed, the time
convolutions in the fourth lines of (\ref{eqn:tractable_u}) and
(\ref{eqn:tractable_v}) can then be viewed as the solution of
inhomogeneous heat equations with very specific inhomogeneous
terms. Properties of solutions to this type of equations are
studied in detail in Section \ref{sect:inhomogeneousheat}
below.

Assuming all results of Section \ref{sect:burgers} and
\ref{sect:inhomogeneousheat}, we now explain how to proceed to
prove that ${\cal F}[{\bf z},{\bf R}^N]$ defines a contraction map.

Obviously, the requirement on $\{g^{\pm}_{n}\}_{n=0}^{N}$ and 
$\{d^{\pm}_n\}_{n=1}^{N}$ is that the first four lines in
(\ref{eqn:tractable_u}) and (\ref{eqn:tractable_v}) satisfy
(\ref{eqn:restestim}). This is achieved in the following way:
\begin{enumerate}
\item \label{item:mass} The first line of
(\ref{eqn:tractable_u}), respectively of (\ref{eqn:tractable_v})
satisfies (\ref{eqn:restestim}) for any
$g_0^{\pm}$ such that the total mass of $g_0^{\pm}$ is equal to
that of $a_0\pm b_0$, provided $a_0\pm b_0$ and $g_0^{\pm}$
satisfy $\|\,x^2(a_0\pm b_0)\|_2<\infty$ and
$\|\,x^2g_0^{\pm}\|_2<\infty$. This fixes the total mass of
$g_0^{\pm}$. Note also that we need the estimate
$\|\,x^2(a_0\pm b_0)\|_2<\infty$. There is no
smallness assumption here, which is to be expected since
generically $\|\,x^2(a(\cdot,t)\pm b(\cdot,t))\|_2$ will grow
as $t\to\infty$. Note on the other hand that
Proposition \ref{prop:weightednorm} shows that
$\|\,x^2(a(\cdot,t)\pm b(\cdot,t))\|_2$ remains finite for all 
$t<\infty$, so requiring $\|\,x^2(a_0\pm b_0)\|_2<\infty$ is
acceptable.

\item We can set the second lines in (\ref{eqn:tractable_u}) and
(\ref{eqn:tractable_v}) equal to zero by picking for $u_0$ and $v_0$
any solution of Burger's equations
\begin{equs}
\partial_t u_0=\partial_x^2u_0+c_{+}\partial_x(u_0)^2~~\mbox{ and }~~
\partial_t v_0=\partial_x^2v_0+c_{-}\partial_x(v_0)^2
\end{equs}
(or of the corresponding heat equations if either $c_{+}$ or
$c_{-}$ happen to be zero). In Proposition \ref{prop:burgers}, we
will prove that there exist unique functions $u_0$ and $v_0$ of
the form given in (\ref{eqn:defu0u1}) that satisfy the conditions
of item \ref{item:mass} above (total mass and decay properties).
This uniquely determines $u_0$ and $v_0$.
\item We can also set the third lines in (\ref{eqn:tractable_u})
and (\ref{eqn:tractable_v}) equal to zero, by picking any
solutions $u_1$ and $v_1$ of linearized Burger's equations
\begin{equs}
\partial_t u_1=\partial_x^2u_1+2c_{+}\partial_x(u_0u_1)
~~\mbox{ and }~~
\partial_t v_1=\partial_x^2v_1+2c_{-}\partial_x(v_0v_1)~.
\label{eqn:unpmvnpm}
\end{equs}
In Proposition \ref{prop:burgers}, we will also prove that there
is a choice of functions $\{g^{\pm}_{n}\}_{n=1}^{N}$ such
that $u_1$ and $v_1$ in (\ref{eqn:defu0u1}) satisfy
(\ref{eqn:unpmvnpm}) for any choice of the coefficients
$\{d^{\pm}_{n}\}_{n=1}^{N}$. Furthermore, in Proposition
\ref{prop:burgers}, we will prove that the choice of
functions can be made in such a way that $g^{\pm}_{n}(x)$ have
Gaussian tails as $x\to\mp\infty$ and algebraic tails as
$x\to\pm\infty$. This actually completely determines
$g^{\pm}_{n}(x)$ up to multiplicative constants (this
last indeterminacy will be removed when the coefficients
$\{d^{\pm}_{n}\}_{n=1}^{N}$ are fixed).
\item \label{item:free}
We then further decompose the terms involving $g_{n}^{\pm}$ in
the fourth lines in (\ref{eqn:tractable_u}) and
(\ref{eqn:tractable_v}) as $g_{n}^{\pm}(x)=f_n(\mp
x)+R^{\pm}_{n}(x)$. The definition and properties of $f_n(x)$
are given in Lemma \ref{lem:alittlelemma}. In particular, in
Proposition \ref{prop:burgers}, we will prove that
$R^{\pm}_{n}(x)$ have zero total mass and Gaussian tails as
$|x|\to\infty$, which implies that $\ed^{\partial_x^2
t}R^{\pm}_{n}$ also satisfy (\ref{eqn:restestim}).
\item\label{item:last} Finally, in Section
\ref{sect:inhomogeneousheat}, we will prove that the time
convolution part of the fourth lines in
(\ref{eqn:tractable_u}) and (\ref{eqn:tractable_v}) can be split
into linear combinations of $\ed^{\partial_x^2 t}f_{n}(\mp x)$
with $n=1\ldots N+1$ plus a remainder that satisfies
(\ref{eqn:restestim}). The coefficients
$\{d^{\pm}_n\}_{n=1}^{N}$ can then be set
recursively by requiring that all the terms with $n=1\ldots N$
coming from the time convolution are canceled by those coming
from item \ref{item:free} above. This can always be done because
the coefficient of $\ed^{\partial_x^2 t}f_{m}(\mp x)$ in the time
convolution part of the fourth lines in (\ref{eqn:tractable_u})
and (\ref{eqn:tractable_v}) depends only on $g_0^{\pm}$ if $m=1$
and on $d^{\pm}_{m-1}$ if $m>1$.
The only term that cannot be set to zero is the last term in the
linear combination (the one with $n=N+1$), which is the one that
`drives' the equations and fixes $\epsilon=2^{-N-2}$.
\end{enumerate}

The procedure outlined in \ref{item:mass}-\ref{item:last} takes
care of the first four lines in (\ref{eqn:tractable_u}) and
(\ref{eqn:tractable_v}). We will then prove in Section
\ref{sect:remainderestimates} that the terms ${\cal
R}_{\{u,v\}}[{\bf z}]$ satisfy (\ref{eqn:restestim}) and that
\begin{equs}
\label{eqn:onRtilde}
\sum_{\alpha=0}^1
\|\D^{\alpha}
\widetilde{\cal R}_{\{u,v\}}
[{\bf z},{\bf R}^{N}]\|_{2,\frac{3}{4}+\frac{\alpha}{2}-\epsilon}
&\leq
C\epsilon_0
\sum_{\alpha=0}^1
\|\D^{\alpha}{\bf R}^{N}\|_{2,\frac{3}{4}+\frac{\alpha}{2}-\epsilon}
+C~,\\
\label{eqn:onRtildeLip}
\sum_{\alpha=0}^1
\|\D^{\alpha}
(\widetilde{\cal R}_{\{u,v\}}[{\bf z},{\bf R}^N_1]
-\widetilde{\cal R}_{\{u,v\}}[{\bf z},{\bf R}^N_2]
)\|_{2,\frac{3}{4}+\frac{\alpha}{2}-\epsilon}&\leq
C\epsilon_0
\sum_{\alpha=0}^1
\|\D^{\alpha}({\bf R}^{N}_1-{\bf R}^{N}_2)\|_{2,\frac{3}{4}
+\frac{\alpha}{2}-\epsilon}~.
\end{equs}
This finally proves that ${\cal F}[{\bf z},{\bf R}^N]$ defines a
contraction map and that the solution of ${\bf R}^N={\cal
F}[{\bf z},{\bf R}^N]$ satisfies (\ref{eqn:restestim}), which
completes the proof of Theorems \ref{thm:maintheorem} and
\ref{thm:asymptoticsrestated}.

\section{Burger's type equations}\label{sect:burgers}

In this section, we consider particular solutions of Burger's type equations
\begin{equs}
\partial_tu_0&=\partial_x^2u_0+\gamma\partial_xu_0^2
\label{eqn:burgers}
\\
\partial_tu_n^{\pm}&=\partial_x^2u_n^{\pm}+2\gamma\partial_x(u_0u_n^{\pm})
\label{eqn:linburgers}
\end{equs}
of the form 
\begin{equs}
{\textstyle
u_0(x,t)=\frac{1}{\sqrt{1+t}}g_0(\frac{x}{\sqrt{1+t}})~~~~~~
\mbox{ and }~~~~~~
u_n^{\pm}(x,t)=\frac{1}{(1+t)^{1-\frac{1}{2^{n+1}}}}
g_n^{\pm}(\frac{x}{\sqrt{1+t}})}~.
\label{eqn:scalingform}
\end{equs}
We will show that for fixed
$\M(u_0)=\intR\hspace{-3mm}{\rm
d}x~u_0(x,t)=\intR\hspace{-3mm}{\rm d}x~g_0(x)$ small enough,
there is a unique choice of $g_0$ and $g_n^{\pm}$ such that
$g_n^{\pm}(x)=f_n(\mp x)+R_n^{\pm}(x)$, where
\begin{equs}
f_{n}(z)=
\int_z^{\infty}
\hspace{-3mm}{\rm d}\xi~
\frac{\xi\ed^{-\frac{\xi^2}{4}}}{(\xi-z)^{1-\frac{1}{2^{n}}}}
\label{eqn:deffn}
\end{equs}
and $R_n^{\pm}$ has zero mean and Gaussian tails as
$|x|\to\infty$. In particular, $g_n^{\pm}(x)$ decays
algebraically as $x\to\pm\infty$, as is apparent from (\ref{eqn:deffn}). 

Before proceeding to our study of (\ref{eqn:burgers}) and
(\ref{eqn:linburgers}), we prove key properties of the functions
$f_n$.
\begin{lemma}\label{lem:alittlelemma}
Fix $1\leq n<\infty$. The function $f_n$ is the unique solution of
\begin{equa}\label{eqn:foneeq}
\partial_z^2f_{n}(z)+{\textstyle\frac{1}{2}}z\partial_zf_{n}(z)
+({\textstyle1-\frac{1}{2^{n+1}}})f_{n}(z)=0~,~~~~\mbox{with}\\
f_{n}(0)=2^{\frac{1}{2^n}}\Gamma({\textstyle\frac{1+2^{-n}}{2}})~~\mbox{ and }~~
\lim_{z\to\infty}z^{-1+\frac{1}{2^n}}\ed^{\frac{z^2}{4}}
f_{n}(z)&<\infty~.
\end{equa}
It satisfies $\intR\hspace{-3mm}{\rm d}z~f_n(z)=0$ and there
exists a constant $C(n)$ such that
\begin{equa}
\sup_{z\in{\bf R}}
\sum_{m=0}^{2}
\rho_{\frac{1}{2^n} -m,1+m-\frac{1}{2^n}}(z)
|\partial_z^m\myl{10}z f_n(z)+2\partial_{z}f_n(z)\myr{10}|&\leq C(n)\\
\sup_{z\in{\bf R}}
\sum_{m=0}^{3}
\rho_{\frac{1}{2^n}-1-m,2+m-\frac{1}{2^n}}(z) |\partial_z^mf_n(z)|
&\leq C(n)~,
\label{eqn:fnestimates}
\end{equa}
where
\begin{equs}
\rho_{p,q}(z)&=\my{\{}{16}
\begin{array}{ll}
(1+z^2)^{\frac{p}{2}}\ed^{\frac{z^2}{4}}&\mbox{ if }z\geq0\\[2mm]
(1+z^2)^{\frac{q}{2}}&\mbox{ if }z\leq0
\end{array}
~.
\end{equs}
\end{lemma}

\begin{proof}
We first note that $f_n$ can be written as
\begin{equs}
f_{n}(z)=
\int_0^{\infty}
\hspace{-3mm}{\rm d}\xi~
\frac{(\xi+z)\ed^{-\frac{(\xi+z)^2}{4}}}{\xi^{1-\frac{1}{2^{n}}}}
=
-2
\int_0^{\infty}
\hspace{-3mm}{\rm d}\xi~
\xi^{\frac{1}{2^{n}}-1}
\partial_{\xi}\myl{10}\ed^{-\frac{(z+\xi)^2}{4}}\myr{10}
~.
\label{eqn:otherfndef}
\end{equs}
This shows that $f_n$ solves (\ref{eqn:foneeq}) since, defining
${\cal L}f\equiv\partial_z^2f+\frac{1}{2}z\partial_zf+
(1-\frac{1}{2^{n+1}})f$, we find
\begin{equs}
{\cal L}f_n(z)=
\int_{0}^{\infty}
\hspace{-3mm}{\rm d}\xi~
\my{[}{12}
\xi^{\frac{1}{2^n}}\partial_{\xi}^2\myl{10}\ed^{-\frac{(z+\xi)^2}{4}}\myr{10}
-{\textstyle\frac{1}{2^{n+1}}}
(-2)\xi^{\frac{1}{2^n}-1}
\partial_{\xi}\myl{10}\ed^{-\frac{(z+\xi)^2}{4}}\myr{10}
\my{]}{12}
=0~.
\end{equs}

Obviously, $f_n(z)$ is finite for all finite $z$, so we only need
to prove that $f_n$ satisfies the correct decay properties as
$|z|\to\infty$ so that (\ref{eqn:fnestimates}) holds. It is
apparent from (\ref{eqn:deffn}) that $f_n$ decays like a
(modified) Gaussian as $z\to\infty$ and algebraically as
$z\to-\infty$. Furthermore, substituting $f(z)=C |z|^{p_1}$ and
$f(z)=C |z|^{p_2}\ed^{-\frac{z^2}{4}}$ into ${\cal L}f=0$ shows
that the only decay rates compatible with ${\cal L}f=0$ are
$p_1=-2+\frac{1}{2^n}$ and $p_2=1-\frac{1}{2^n}$.

We now complete the proof of the decay estimates
(\ref{eqn:fnestimates}). Let
$F_{n,m}(\xi,z)=\partial_z^m((\xi+z)\ed^{-\frac{(\xi+z)^2}
{4}})$ and
$G_{n,m}(\xi,z)=\partial_z^m(zF_n(\xi,z)+2\partial_zF_n(\xi,z))$.

We first consider the case $z>0$ and note that $F_{n,m}$ and
$G_{n,m}$ satisfy
\begin{equs}[2]
|F_{n,m}(\xi,z)|&\leq |F_{n,m}(0,z)|&~~~\mbox{ and }~~~
|G_{n,m}(\xi,z)|&\leq |G_{n,m}(0,z)|
\end{equs}
for all $\xi\geq0$ if $z\geq z_0$ for some $z_0$ large enough. We
thus get, e.g.
\begin{equs}
|f_{n}(z)|&=
\my{|}{14}
\int_0^{\infty}
\hspace{-3mm}{\rm d}\xi~
F_{n,0}(\xi,z)\xi^{\frac{1}{2^n}-1}
\my{|}{14}\leq
|F_{n,0}(0,z)|
\int_0^{z^{-1}}
\hspace{-3mm}{\rm d}\xi~
\xi^{\frac{1}{2^n}-1}+
z^{1-\frac{1}{2^n}}
\int_{z^{-1}}^{\infty}
\hspace{-3mm}{\rm d}\xi~
|F_{n,0}(\xi,z)|
\leq C z^{1-\frac{1}{2^n}}\ed^{-\frac{z^2}{4}}~.
\end{equs}
The estimates on $|\partial_z^m(zf_n(z)+2\partial_{z}f_n(z))|$
and $|\partial_z^{1+m} f_n(z)|$ when $z>0$ and $m\geq1$ can be done in
exactly the same way; hence we omit the details.             

We now consider the case $z<0$ and note that $F_{n,m}$ and
$G_{n,m}$ satisfy
\begin{equs}[2]
|F_{n,m}(\xi,z)|&\leq |F_{n,m}(-{\textstyle \frac{z}{2}},z)|~~
&\mbox{ and }~~
|G_{n,m}(\xi,z)|&\leq |G_{n,m}(-{\textstyle \frac{z}{2}},z)|
\end{equs}
for all $0\leq\xi\leq-\frac{z}{2}$ if $z\leq-z_0$ for some $z_0$ large enough.
We thus find (integrating by parts in the second integral below)
\begin{equs}
|f_{n}(z)|&=
\my{|}{14}
\int_0^{\infty}
\hspace{-3mm}{\rm d}\xi~
F_{n,0}(\xi,z)\xi^{\frac{1}{2^n}-1}
\my{|}{14}\leq
|F_{n,0}(-{\textstyle \frac{z}{2}},z)|
\int_0^{-\frac{z}{2}}
\hspace{-3mm}{\rm d}\xi~
\xi^{\frac{1}{2^n}-1}+
\my{|}{14}
\int_{-\frac{z}{2}}^{\infty}
\hspace{-3mm}{\rm d}\xi~
F_{n,0}(\xi,z)
\xi^{\frac{1}{2^n}-1}\my{|}{14}
\\&
\leq C
|z|^{\frac{1}{2^n}-1}\ed^{-\frac{z^2}{16}}
+
2\myl{10}
1-
{\textstyle\frac{1}{2^n}}
\myr{10}
\int_{-\frac{z}{2}}^{\infty}
\hspace{-3mm}{\rm d}\xi~
\ed^{-\frac{(\xi+z)^2}{4}}\xi^{\frac{1}{2^n}-2}
\leq C|z|^{\frac{1}{2^n}-2}~.
\end{equs}
Since the remaining estimates can again be done in exactly the
same way, we omit the details. It only remains to show that
$f_n(z)$ has zero total mass. This follows from
\begin{equs}
\int_{-\infty}^{\infty}
\hspace{-3mm}{\rm d}z~
f_n(z)=
(
{\textstyle\frac{1}{2}-\frac{1}{2^{n+1}}}
)^{-1}\int_{-\infty}^{\infty}
\hspace{-3mm}{\rm d}z~
{\cal L}f_n(z)
=0
~,
\end{equs}
since $\partial_z^2f_n,z\partial_zf_n$ and $f_n$ are all integrable
over ${\bf R}$.
\end{proof}

\begin{remark}
Using the representation (\ref{eqn:otherfndef}), splitting the
integration interval into $[0,2^{-\frac{n}{2}})$ and
$[2^{-\frac{n}{2}},\infty)$, integrating by parts and letting
$n\to\infty$, one can prove that
\begin{equs}
\lim_{n\to\infty}}2^{-n}f_n(z)=z\ed^{-\frac{z^2}{4}~,
\end{equs}
which shows that the constant $C(n)$ in (\ref{eqn:fnestimates})
grows at most like $2^n$.
\end{remark}

We can now study in detail the solutions of (\ref{eqn:burgers})
and (\ref{eqn:linburgers}) that are of the form (\ref{eqn:scalingform}):
\begin{proposition}
\label{prop:burgers}
Fix $1\leq n<\infty$.
For all $\alpha,\gamma\in{\bf R}$ with $|\alpha\gamma|$
small enough, there exist unique functions $u_0$ and $u_n^{\pm}$
of the form (\ref{eqn:scalingform}) that solve
(\ref{eqn:burgers}) and (\ref{eqn:linburgers}), with $g_0$
satisfying
\begin{equs}
\int_{-\infty}^{\infty}\hspace{-3mm}{\rm d}z~g_0(z)=\alpha~,~~~
\sum_{m=0}^{3}
\frac{\ed^{\frac{z^2}{4}}}{
(\sqrt{1+z^2})^{m}}|\partial_z^mg_0(z)|
\leq C |\alpha|
\end{equs}
and with $g_n^{\pm}(z)=f_{n}(\mp z)+R_n^{\pm}(z)$, where
$R_n^{\pm}$ satisfy
\begin{equs}
\int_{-\infty}^{\infty}
\hspace{-3mm}{\rm d}z~R_n^{\pm}(z)=0~~\mbox{ and }~~
{\displaystyle\sup_{z\in{\bf R}}
\sum_{m=0}^{3}}
\frac{\ed^{\frac{z^2}{4}}}{
(\sqrt{1+z^2})^{1+m-\frac{1}{2^n}}}|\partial_z^mR_n^{\pm}(z)|
\leq C |\alpha\gamma|~.
\end{equs}
\end{proposition}

\begin{proof}
The (unique) solution of (\ref{eqn:burgers}) of the form
$u_0(x,t)=\frac{1}{\sqrt{1+t}}g_0(\frac{x}{\sqrt{1+t}})$ 
satisfying $\intR\hspace{-3mm}{\rm d}z~g_0(z)=\alpha$
is given by 
\begin{equs}
g_0(z)=\frac{
\tanh(\frac{\alpha\gamma}{2})\ed^{-\frac{z^2}{4}}
}{
\gamma\sqrt{\pi}
(1+\tanh(\frac{\alpha\gamma}{2})\erf(\frac{z}{2}))
}~.
\end{equs}
In particular, we have
\begin{equs}
\sum_{m=0}^{3}
\frac{\ed^{\frac{z^2}{4}}}{
(\sqrt{1+z^2})^{m}}|\partial_z^mg_0(z)|
\leq C |\alpha|~.
\label{eqn:f0bound}
\end{equs}
We next note that substituting (\ref{eqn:scalingform}) into
(\ref{eqn:linburgers}) gives
\begin{equs}
0&=
{\textstyle\partial_z^2g_n^{\pm}(z)+\frac{1}{2}z\partial_zg_n^{\pm}(z)
+(1-\frac{1}{2^{n+1}})g_n^{\pm}(z)}+
2\gamma\partial_z(g_0(z)g_n^{\pm}(z))\\[2mm]
&\equiv{\cal L}g_n^{\pm}(z)+2\gamma\partial_z(u_0(z)g_n^{\pm}(z))~.
\label{eqn:fpmeq}
\end{equs}
We formally have (using integration by parts)
\begin{equs}
\int_{-\infty}^{\infty}
\hspace{-3mm}{\rm d}z~
g_n^{\pm}(z)=
(
{\textstyle\frac{1}{2}-\frac{1}{2^{n+1}}}
)^{-1}\int_{-\infty}^{\infty}
\hspace{-3mm}{\rm d}z~
{\cal L}g_n^{\pm}(z)
+2\gamma\partial_z(u_0(z)g_n^{\pm}(z))
=
0
~,
\label{eqn:formal}
\end{equs}
which shows that $g_n^{\pm}$ have zero total mass, {\em provided
the formal manipulations above are justified}, i.e. provided
$g_n^{\pm}$ and its derivatives decay fast enough so that the
integrals are convergent.

As is easily seen, $f_n(z)$ and $f_n(-z)$ are two linearly
independent solutions of ${\cal L}f=0$, whose general solution
can thus be written as $c_1f_n(z)+c_2f_n(-z)$. Using the variation of
constants formula, we get that the solution of (\ref{eqn:fpmeq})
satisfies the integral equation
\begin{equs}
g_n^{\pm}(z)
=f_n(z)\myl{14}
c_{1}^{\pm}+2\gamma\int_0^{z}
\hspace{-2mm}{\rm d}\xi~
{\textstyle\frac{
f_n(-\xi)\partial_{\xi}(g_0(\xi)g_n^{\pm}(\xi))
}{
W(\xi)}
}
\myr{14}+
f_n(-z)\myl{14}
c_{2}^{\pm}-2\gamma\int_0^{z}
\hspace{-2mm}{\rm d}\xi~
{\textstyle\frac{
f_n(\xi)\partial_{\xi}(g_0(\xi)g_n^{\pm}(\xi))
}{
W(\xi)
}}
\myr{14}~,
\end{equs}
where the Wronskian $W(z)$ is given by
$W(z)=f_n(z)\partial_zf_n(-z)-f_n(-z)\partial_zf_n(z)$ and
$c_{1}^{\pm}$ and $c_{2}^{\pm}$ are free parameters. Note that
$W(z)$ satisfies $\partial_zW(z)=-\frac{z}{2}W(z)$ and hence
$W(z)=W(0)\ed^{-\frac{z^2}{4}}$ for some $W(0)\neq0$. We now set
$c_{1}^{\pm}$ and $c_{2}^{\pm}$ in such a way that (after integration by
parts), we have
\begin{equs}
g_n^{\pm}(z)&=f_n(\mp z)
+R[g_n^{\pm}](z)~,
\label{eqn:contractforpm}
\\
R[g_n^{\pm}](z)&=
{\textstyle\frac{\gamma}{W(0)}}
f_n(z)
\int_{-\infty}^{z}
\hspace{-3mm}{\rm d}\xi~
\ed^{\frac{\xi^2}{4}}
(\xi f_n(-\xi)+2\partial_{\xi}f_n(-\xi))
g_0(\xi)g_n^{\pm}(\xi)\\
&\phantom{=}~+
{\textstyle\frac{\gamma}{W(0)}}
f_n(-z)
\int_{z}^{\infty}
\hspace{-3mm}{\rm d}\xi~
\ed^{\frac{\xi^2}{4}}
(\xi f_n(\xi)+2\partial_{\xi}f_n(\xi))
g_0(\xi)g_n^{\pm}(\xi)~.
\end{equs}
Using Lemma \ref{lem:alittlelemma} and (\ref{eqn:f0bound}), it is
then easy to show that for $|\alpha\gamma|$ small enough,
(\ref{eqn:contractforpm}) defines a contraction map in the norm
\begin{equs}
|f|_{2-\frac{1}{2^n}}\equiv
\sup_{z\in{\bf R}}(\sqrt{1+z^2})^{2-\frac{1}{2^{n}}}|f(z)|~.
\end{equs}
Namely, we have the improved decay rates
\begin{equs}
\sup_{z\in{\bf R}}
\sum_{m=0}^1
\frac{\ed^{\frac{z^2}{4}}}{
(\sqrt{1+z^2})^{1+m-\frac{1}{2^n}}}|\partial_z^mR[g_n^{\pm}](z)|
&\leq
C|\alpha\gamma|~
|g_n^{\pm}|_{2-\frac{1}{2^n}}~.
\end{equs}
This shows that (\ref{eqn:contractforpm}) has a (locally) unique
solution among functions with 
$|f|_{2-\frac{1}{2^n}}\leq c_0$ if $|\alpha\gamma|$ is
small enough. In particular, there holds
\begin{equs}
\sup_{z\in{\bf R}}
\sum_{m=0}^1
\frac{\ed^{\frac{z^2}{4}}}{
(\sqrt{1+z^2})^{1+m-\frac{1}{2^n}}}|\partial_z^mR[g_n^{\pm}](z)|
&\leq
C|\alpha\gamma|
~,
\end{equs}
from which we deduce, using again (\ref{eqn:contractforpm}) and
Lemma \ref{lem:alittlelemma}, that $|\D g_n^{\pm}|_{3-\frac{1}{2^n}}\leq c_1$
and thus
\begin{equs}
\sup_{z\in{\bf R}}
\frac{\ed^{\frac{z^2}{4}}}{
(\sqrt{1+z^2})^{3-\frac{1}{2^n}}}|\partial_z^2R[g_n^{\pm}](z)|
&\leq
C|\alpha\gamma|
~.
\end{equs}
Iterating this procedure shows that $|\D^mg_n^{\pm}|_{2+m-\frac{1}{2^n}}\leq c_m$ and that
\begin{equs}
\sup_{z\in{\bf R}}
\sum_{m=0}^3
\frac{\ed^{\frac{z^2}{4}}}{
(\sqrt{1+z^2})^{1+m-\frac{1}{2^n}}}|\partial_z^mR[g_n^{\pm}](z)|
&\leq
C|\alpha\gamma|
\end{equs}
as claimed. In turn, this proves that the formal manipulations in
(\ref{eqn:formal}) are justified, so that the functions
$g_n^{\pm}(z)$ have zero total mass, which shows that
the remainders $R[g_n^{\pm}](z)$ have zero total mass as claimed
since $R[g_n^{\pm}](z)=g_n^{\pm}(z)-f_n(\pm z)$ and
since both $g_n^{\pm}(z)$ and $f_n(z)$ have zero total mass.
\end{proof}

\section{Inhomogeneous heat equations}\label{sect:inhomogeneousheat}

In this section, we consider solutions of inhomogeneous heat
equations of the form
\begin{equs}
\partial_t u=\partial_x^2u+\partial_x
\myl{13}
(1+t)^{\frac{1}{2^n}-\frac{3}{2}}
f\myl{10}
{\textstyle
\frac{x-2\sigma t}{\sqrt{1+t}}
}
\myr{10}
\myr{13}~,~~~~u(x,0)=0~,
\label{eqn:inhomoheateqgeneral}
\end{equs}
where $f$ is a regular function having Gaussian decay at
infinity. Solutions of (\ref{eqn:inhomoheateqgeneral}) satisfy
\begin{theorem}
Let $1\leq n<\infty$, $\sigma=\pm1$,
$\Xi(x)=\ed^{\frac{x^2}{8}}$, $\M(f)=
\intR\hspace{-3mm}{\rm d}z~f(z)$ and
\begin{equs}
u_{n}(x,t)={\textstyle\frac{\sigma}{(1+t)^{1-\frac{1}{2^{n+1}}}}
\frac{2^{-1-\frac{1}{2^n}}}{\sqrt{4\pi}}
f_n(\frac{-\sigma x}{\sqrt{1+t}})}~~~~\mbox{ with }~~~~
f_n(z)=
\int_z^{\infty}
\hspace{-3mm}{\rm d}\xi~
\frac{\xi\ed^{-\frac{\xi^2}{4}}}{(\xi-z)^{1-\frac{1}{2^n}}}~.
\label{eqn:defun}
\end{equs}
The solution $u$ of (\ref{eqn:inhomoheateqgeneral}) satisfies
\begin{equs}
\|u-\M(f)\,u_{n}\|_{2,\frac{3}{4}^{\star}}
+
\|\D(u-\M(f)\,u_{n})\|_{2,\frac{5}{4}^{\star}}
\leq C
\sum_{m=0}^2\|\Xi\D^m f\|_{\infty}
~,
\label{eqn:estiinhomoheat}
\end{equs}
for all $f$ such that the r.h.s. of (\ref{eqn:estiinhomoheat}) is
finite.
\end{theorem}

\begin{remark}
Note that while $u\to \M(f)u_{n}$ as $t\to\infty$ in the Sobolev
norm (\ref{eqn:estiinhomoheat}), it does not do so in spatially
weighted norms such as $\L^2({\bf R},x^2\d x)$, as $u_{n}$ has
infinite spatial moments for all times, while all moments of $u$
are bounded for finite time.
\end{remark}

\begin{proof}
We first define
\begin{equs}
F(\xi)=
\int_{-\infty}^{\xi}
\hspace{-3mm}{\rm d}z~\myl{10}
f(z)-
\M(f)\,{\textstyle
\frac{\ed^{-\frac{z^2}{4}}}{\sqrt{4\pi}}}
\myr{10}~~~~\mbox{ with }~~~~
\M(f)=\int_{-\infty}^{\infty}
\hspace{-3mm}{\rm d}z~f(z)
\label{eqn:defF}
\end{equs}
and note that $F$ satisfies
\begin{equs}
 \|\D^3F\|_1
+\sum_{m=0}^{2}\|\rho\D^m F\|_1
+\sum_{m=1}^{2}\|\D^m F\|_2
\leq
C
\sum_{m=0}^2\|\Xi\D^m f\|_{\infty}
~,
\label{eqn:FF}
\end{equs}
where $\rho(x)=\sqrt{1+x^2}$. Namely, we first note that $\|\rho
F\|_1\leq\|\hat{F}\|_2+\|\hat{F}''\|_2$ and
$\hat{F}(k)=(ik)^{-1}(\hat{f}(k)-\hat{f}(0)\ed^{-k^2})$. Then,
since $\|\Xi f\|_{\infty}<\infty$ implies that $\hat{f}$ is
analytic, $\hat{F}$ is regular near $k=0$. The proof of
(\ref{eqn:FF}) now follows from elementary arguments.

We finally note that it follows from (\ref{eqn:defF}) that
\begin{equs}
(1+t)^{\frac{1}{2^n}-\frac{3}{2}}
f\myl{10}
{\textstyle
\frac{x-2\sigma t}{\sqrt{1+t}}
}
\myr{10}
=
\M(f)~
\underbrace{\frac{(1+t)^{\frac{1}{2^n}-\frac{3}{2}}}{\sqrt{4\pi}}
\ed^{-\frac{(x-2\sigma t)^2}{4(1+t)}}}_{\equiv A(x,t)}+
\underbrace{
(1+t)^{\frac{1}{2^n}-1}
\partial_x
F\myl{10}
{\textstyle
\frac{x-2\sigma t}{\sqrt{1+t}}
}
\myr{10}}_{\equiv \partial_x B(x,t)}~.
\label{eqn:underbraces}
\end{equs}
The proof of (\ref{eqn:estiinhomoheat}) is then completed by
considering separately the solutions of heat equations with
inhomogeneous terms given by $\partial_x A(x,t)$ and
$\partial_x^2 B(x,t)$. This is done in Propositions
\ref{prop:Gaussian} and \ref{prop:secondder} below.
\end{proof}

\begin{proposition}
\label{prop:Gaussian}
Let $\sigma=\pm1$, $1\leq n<\infty$, and let $u_n$ be defined as in
(\ref{eqn:defun}). The solution $u$ of
\begin{equs}
\partial_t u=\partial_x^2u+\partial_xA
~,~~~~u(x,0)=0~,
\label{eqn:inhomoheateq}
\end{equs}
with $A$ defined in (\ref{eqn:underbraces}) satisfies
\begin{equs}
\|u-u_{n}\|_{2,\frac{3}{4}}
+
\|\D(u-u_{n})\|_{2,\frac{5}{4}}
\leq C~.
\label{eqn:estiinhomoheatrap}
\end{equs}
\end{proposition}

\begin{proof}
The solution of (\ref{eqn:inhomoheateq}) is given by
\begin{equs}
u(x,t)=
\partial_x
\int_0^t
\hspace{-2mm}
{\rm d}s
\int_{-\infty}^{\infty}
\hspace{-3mm}
{\rm d}y
\frac{\ed^{-\frac{(x-y)^2}{4(t-s)}}}{\sqrt{4\pi(t-s)}}
\frac{\ed^{-\frac{(y-2\sigma s)^2}{4(1+s)}}}{\sqrt{4\pi}(1+s)^{\frac{3}{2}-\frac{1}{2^{n}}}}~.
\label{eqn:inhomoheat}
\end{equs}
To motivate our result, we note that performing the
$y$-integration and changing variables from 
$s$ to $\xi\equiv\frac{2s-\sigma x}{\sqrt{1+t}}$ in
(\ref{eqn:inhomoheat}) leads to
\begin{equs}
\lim_{t\to\infty}
(1+t)^{1-\frac{1}{2^{n+1}}}u(-\sigma z\sqrt{1+t},t)=
\lim_{t\to\infty}
{\textstyle\frac{\sigma 2^{-1-{\frac{1}{2^n}}}}{\sqrt{4\pi}}}
\int_{z}^{\frac{2t}{\sqrt{1+t}}+z}
\hspace{-3mm}{\rm d}\xi~
{\textstyle\frac{
\xi\ed^{-\frac{\xi^2}{4}}
}{
(\xi-z+\frac{2}{\sqrt{1+t}})^{1-\frac{1}{2^n}}
}}                                                              
={\textstyle\frac{\sigma 2^{-1-{\frac{1}{2^n}}}}{\sqrt{4\pi}}}
f_n(z)~.
\end{equs}
More formally, taking the Fourier transform of
(\ref{eqn:inhomoheat}) gives
\begin{equs}
\hat{u}(k,t)&=ik\ed^{-k^2(1+t)}
\int_0^{t}
\hspace{-2mm}{\rm d}s~
\frac{\ed^{2ik\sigma s}}{(1+s)^{1-\frac{1}{2^n}}}~.
\end{equs}
We now use that
\begin{equs}
\my{|}{14}
\int_0^{t}
\hspace{-2mm}{\rm d}s~
\frac{\ed^{2ik\sigma s}}{(1+s)^{1-\frac{1}{2^n}}}
-\int_0^{t}
\hspace{-2mm}{\rm d}s~
\frac{\ed^{2ik\sigma s}}{s^{1-\frac{1}{2^n}}}
\my{|}{14}&\leq C(n)~,\\
\int_0^{t}
\hspace{-2mm}{\rm d}s~
\frac{\ed^{2ik\sigma s}}{s^{1-\frac{1}{2^n}}}&=
|k|^{-\frac{1}{2^n}}
\myl{10}
\theta(\sigma k)J_n(|k|t)+\theta(- \sigma k)\overline{J_n(|k|t)}
\myr{10}~,
\end{equs}
where $\theta(k)$ is the Heaviside step function and we defined
\begin{equs}
J_n(z)=
\int_0^{z}
\hspace{-2mm}{\rm d}s~
\frac{\ed^{2is}}{s^{1-\frac{1}{2^n}}}
\end{equs}
for $z\geq0$. This function satisfies
\begin{equs}
\sup_{z\geq0}
z^{1-\frac{1}{2^n}}
|J_n(z)-J_{n,\infty}|\leq\frac{1}{2}~~\mbox{ for }~~
J_{n,\infty}=\lim_{z\to\infty}J_n(z)~.
\end{equs}
Now define
\begin{equs}
\widehat{u_{n}}(k,t)=
ik\ed^{-k^2(1+t)}
|k|^{-\frac{1}{2^n}}
\myl{10}
\theta(\sigma k)J_{n,\infty}+\theta(-\sigma k)\overline{J_{n,\infty}}
\myr{10}~.
\label{eqn:widehat}
\end{equs}
We have
\begin{equs}
|\hat{u}(k,t)-\widehat{u_{n}}(k,t)|\leq 
(C(n)|k|+t^{-1+\frac{1}{2^n}})
\ed^{-k^2(1+t)}
\leq 
(C(n)|k|+t^{-\frac{1}{2}})
\ed^{-k^2(1+t)}~,
\label{eqn:hammer}
\end{equs}
from which (\ref{eqn:estiinhomoheatrap}) follows by direct integration.
We complete the proof by showing that the inverse Fourier transform
of the function $\widehat{u_{n}}(k,t)$ defined in
(\ref{eqn:widehat}) satisfies
\begin{equs}
u_{n}(x,t)={\textstyle\frac{\sigma}{(1+t)^{1-\frac{1}{2^{n+1}}}}
\frac{2^{-1-\frac{1}{2^n}}}{\sqrt{4\pi}}
f_n(\frac{-\sigma x}{\sqrt{1+t}})}~~\mbox{ for }~~
f_n(z)=\int_z^{\infty}
\hspace{-3mm}{\rm d}\xi~
\frac{\xi\ed^{-\frac{\xi^2}{4}}}{(\xi-z)^{1-\frac{1}{2^n}}}~.
\label{eqn:tocheck}
\end{equs}
This follows easily from the fact that
\begin{equs}
\widehat{u_{n}}(k,t)=(1+t)^{-\frac{1}{2}+\frac{1}{2^{n+1}}}
\widehat{u_{n}}(k\sqrt{1+t},0)~,
\end{equs}
and that, since
\begin{equs}
f_n(z)=
\int_0^{\infty}
\hspace{-3mm}{\rm d}\xi~
\frac{(z+\xi)\ed^{-\frac{(z+\xi)^2}{4}}}{\xi^{1-\frac{1}{2^n}}}~,
\end{equs}
we get
\begin{equs}
{\textstyle\frac{\sigma2^{-1-\frac{1}{2^n}}}{\sqrt{4\pi}}}
\widehat{f_n}(-\sigma k)&=2^{-\frac{1}{2^n}}i k\ed^{-k^2}
\int_0^{\infty}
\hspace{-3mm}{\rm d}\xi~
\frac{\ed^{ik\sigma \xi}}{\xi^{1-\frac{1}{2^n}}}
=ik\ed^{-k^2}|k|^{-\frac{1}{2^n}}
\int_0^{\infty}
\hspace{-3mm}{\rm d}\xi~
\frac{\ed^{2i\sign(k\sigma)\xi}}{\xi^{1-\frac{1}{2^n}}}\\
&=
ik\ed^{-k^2}|k|^{-\frac{1}{2^n}}
\myl{10}
\theta(k\sigma)J_{n,\infty}+
\theta(-k\sigma)\overline{J_{n,\infty}}
\myr{10}
=\widehat{u_{n}}(k,0)
\end{equs}
as claimed.
\end{proof}

\begin{proposition}
\label{prop:secondder}
Let $\sigma=\pm1$, $1\leq n<\infty$ and $\rho(x)=\sqrt{1+x^2}$.
The solution $u$ of
\begin{equs}
\partial_t u=\partial_x^2u+\partial_x^2 B
~,~~~~u(x,0)=0~,
\label{eqn:inhomoheateqF}
\end{equs}
with $B$ defined in (\ref{eqn:underbraces}) satisfies
\begin{equs}
 \|u\|_{2,\frac{3}{4}^{\star}}
+\|\D u\|_{2,\frac{5}{4}^{\star}}
\leq C~
\myl{13}
 \|\D^3F\|_1
+\sum_{m=0}^{2}\|\rho\D^m F\|_1
+\sum_{m=1}^{2}\|\D^m F\|_2
\myr{13}
\label{eqn:estiinhomoheatF}
\end{equs}
for all $F$ for which the r.h.s. of (\ref{eqn:estiinhomoheatF})
is finite.
\end{proposition}

\begin{proof}
We first note that the Fourier transform of $u$ is given by
\begin{equs}
\hat{u}(k,t)&=
-k^2
\int_0^t
\hspace{-2mm}{\rm d}s~
\ed^{-k^2(t-s)-2ik\sigma s}
\hat{F}(k\sqrt{1+s})(1+s)^{\frac{1}{2^n}-\frac{1}{2}}~,
\end{equs}
which implies
\begin{equs}
\|(1-\Q)u\|_{2,\frac{3}{4}}
+
\|(1-\Q)\D u\|_{2,\frac{5}{4}}
\leq C\myl{10}\|\D F\|_2+\|\D^2 F\|_2\myr{10}
\sup_{0\leq t\leq1}
\int_0^t
\hspace{-2mm}
\frac{{\rm d}s}{\sqrt{t-s}}
~.
\end{equs}
Here $\Q$ is again defined as the characteristic function for
$t\geq1$. Next, integrating by parts, we find
\begin{equs}
\hat{u}(k,t)&=
 \frac{ik\hat{F}(k)\ed^{-k^2t}}{2\sigma}
-\frac{ik\hat{F}(k\sqrt{1+t})\ed^{-2ik\sigma t}}{2\sigma
(1+t)^{\frac{1}{2}-\frac{1}{2^n}}
}
+\hat{N}(k,t)\\
\mbox{where }~~~
\hat{N}(k,t)&=
\frac{ik}{2\sigma}
\int_0^t
\hspace{-2mm}{\rm d}s~
\ed^{-k^2(t-s)-2ik\sigma s}
\myl{10}k^2+\partial_s\myr{10}~
\myl{14}
\frac{\hat{F}(k\sqrt{1+s})}{(1+s)^{\frac{1}{2}-\frac{1}{2^n}}}
\myr{14}~.
\end{equs}
We then note that
\begin{equs}
\|u-N\|_{2,\frac{3}{4}}
+
\|\D(u-N)\|_{2,\frac{5}{4}}
\leq C\myl{10}\|F\|_1+\|\D F\|_2+\|\D^2 F\|_2\myr{10}~, 
\end{equs}
and that, defining $\hat{G}(k)=\frac{1}{2}\partial_k\hat{F}(k)$,
we have $\hat{N}(k,t)=\hat{N}_0(k,t)+\hat{N}_1(k,t)+\hat{N}_2(k,t)$, where
\begin{equs}
\hat{N}_0(k,t)&=
\frac{ik^3}{2\sigma}
\int_0^t
\hspace{-2mm}{\rm d}s~
\ed^{-k^2(t-s)-2ik\sigma s}
~
\myl{14}
\frac{\hat{F}(k\sqrt{1+s})}{(1+s)^{\frac{1}{2}-\frac{1}{2^n}}}
\myr{14}~,\\
\hat{N}_1(k,t)&=
\frac{ik^2}{2\sigma}
\int_0^t
\hspace{-2mm}{\rm d}s~
\ed^{-k^2(t-s)-2ik\sigma s}
~
\myl{14}
\frac{\hat{G}(k\sqrt{1+s})}{(1+s)^{1-\frac{1}{2^n}}}
\myr{14}~,\\
\hat{N}_2(k,t)&=
\frac{ik}{2\sigma}
\myl{10}
{\textstyle
\frac{1}{2^n}-\frac{1}{2}
}
\myr{10}
\int_0^t
\hspace{-2mm}{\rm d}s~
\ed^{-k^2(t-s)-2ik\sigma s}
~
\myl{14}
\frac{\hat{F}(k\sqrt{1+s})}{(1+s)^{\frac{3}{2}-\frac{1}{2^n}}}
\myr{14}~.
\end{equs}
The procedure is now similar to that outlined in the proof of
Theorem \ref{thm:cauchy}: split the integration intervals
into $[0,\frac{t}{2}]$ and $[\frac{t}{2},t]$ and distribute the
derivatives ($k$-factors) either on the functions $F$ and $G$, or on the
Gaussian. Introducing the notation
\begin{equs}
\Bone{p_1,q_1}{p_2,q_2}
\equiv
\int_0^{\frac{t}{2}}
\hspace{-3mm}{\rm d}s~
\frac{
(1+s)^{-q_1}
}{
(t-s)^{p_1}
}+
\int_{\frac{t}{2}}^t
\hspace{-2mm}{\rm d}s~
\frac{
(1+s)^{-q_2}
}{
(t-s)^{p_2}
}
~,
\label{eqn:defB}
\end{equs}
we then find that 
\begin{equs}
\|\Q\D^{\alpha}N_0\|_{2,\frac{3}{4}+\frac{\alpha}{2}}
&\leq
C(\|F\|_1+\|\D^{2+\alpha} F\|_1)
\sup_{t\geq1}
t^{\frac{3}{4}+\frac{\alpha}{2}}
~\Bone{\frac{7}{4}+\frac{\alpha}{2},0}{\frac{3}{4},1+\frac{\alpha}{2}}
~,\\
\|\Q\D^{\alpha}N_1\|_{2,\frac{3}{4}+\frac{\alpha}{2}}
&\leq
C(\|G\|_1+\|\D^{1+\alpha} G\|_1)
\sup_{t\geq1}
t^{\frac{3}{4}+\frac{\alpha}{2}}
~\Bone{\frac{5}{4}+\frac{\alpha}{2},\frac{1}{2}}
   {\frac{3}{4},1+\frac{\alpha}{2}}~,\\
\|\Q\D^{\alpha}N_2\|_{2,\frac{3}{4}+\frac{\alpha}{2}^{\star}}
&\leq
C(\|F\|_1+\|\D^{\alpha} F\|_1)
\sup_{t\geq1}
\frac{t^{\frac{3}{4}+\frac{\alpha}{2}}}{\ln(2+t)}
~\Bone{\frac{3}{4}+\frac{\alpha}{2},1}
   {\frac{3}{4},1+\frac{\alpha}{2}}
\end{equs}
for $\alpha=0,1$. The proof is completed by a straightforward
application of Lemma \ref{lem:onB} below, where we consider
generalizations of the function $B_1$ in (\ref{eqn:defB}) (see
Definition \ref{def:defB} below), since those will occur later on
in Sections \ref{sect:cauchyproof} and
\ref{sect:remainderestimates}.
\end{proof}

\section{Proof of Theorem \ref{thm:cauchy}, continued}
\label{sect:cauchyproof}
In view of the estimates (\ref{eqn:estimateonDkernel}) and
(\ref{eqn:H}) on $\ed^{\L t}$ and $h$, respectively, the
estimates needed to conclude the proof of Theorem
\ref{thm:cauchy} will naturally involve the functions $B_0$ and
$B$ which are defined as follows:
\begin{definition}\label{def:defB}
We define
\begin{equs}
B_0[q](t)&=
\int_0^{t}
\hspace{-2mm}{\rm d}s
\frac{\ed^{-\frac{t-s}{8}}}{\sqrt{t-s}(1+s)^q}~,\\
\B{p_1,q_1,r_1}
  {p_2,q_2,r_2,r_3}
&=
\int_0^{\frac{t}{2}}
\hspace{-3mm}{\rm d}s~
\frac{
(1+s)^{-q_1}
}{
(t-s)^{p_1}(t-1+s)^{r_1}
}+
\int_{\frac{t}{2}}^t
\hspace{-2mm}{\rm d}s~
\frac{
(1+s)^{-q_2}
\ln(2+s)^{r_3}
}{
(t-s)^{p_2}(t-1+s)^{r_2}
}~.
\label{eqn:defgeneralB}
\end{equs} 
\end{definition}
These functions satisfy the following estimates:
\begin{lemma}
\label{lem:onB}
Let $0\leq p_2<1$, $0\leq r_2\leq1-p_2$, $p_1,q_1,q_2,r_1\geq0$ and $r_3\in\{0,1\}$. There
exists a constant $C$ such that for all
$t\geq0$ there holds
\begin{equs}
B_0[q_1](t)&\leq
C(1+t)^{-q_1}~,\\
\B{p_1,q_1,r_1}
  {p_2,q_2,r_2,r_3}
&\leq
C~
\ln(2+t)^{\alpha}
\my{\{}{18}
\begin{array}{ll}
\frac{1}{(1+t)^{\beta}}
&~~\mbox{ if }~~0\leq p_1\leq 1
\\
\frac{1}{t^{p_1-1}~(1+t)^{\beta-p_1+1}}
&~~\mbox{ if }~~p_1>1
\end{array}
~,
\label{eqn:generalB}
\end{equs}
where $\beta=\min(p_1+\min(q_1-1,0)+r_1,p_2+q_2+r_2-1)$, 
$\alpha=\max(\delta_{q_1,1},\delta_{p_2+r_2,1}+r_3)$ and
$\delta_{i,j}$ is the Kronecker delta. Furthermore, since
\begin{equs}
\Bone{p_1,q_1}{p_2,q_2}=
\B{p_1,q_1,0}
  {p_2,q_2,0,0}~,
\end{equs}
the estimate in (\ref{eqn:generalB}) applies for
$B_1$ as well.
\end{lemma}

\begin{proof}
The proof follows immediately from
\begin{equs}
B_0[q_1](t)&\leq
\ed^{-\frac{t}{16}}
\int_0^{\frac{t}{2}}
\hspace{-3mm}
\frac{{\rm d}s}{\sqrt{t-s}}
+
\frac{
1
}{
(\frac{t}{2}+1)^{q_1}
}
\int_0^{\frac{t}{2}}
\hspace{-3mm}{\rm d}s~
\frac{\ed^{-\frac{s}{8}}}{\sqrt{s}}
~,\\
\B{p_1,q_1,r_1}
  {p_2,q_2,r_2,r_3}
&\leq
\frac{
1
}{
(\frac{t}{2})^{p_1}(\frac{t}{2}+1)^{r_1}
}
\int_0^{\frac{t}{2}}
\hspace{-3mm}
\frac{{\rm d}s}{(1+s)^{q_1}}
+
\frac{
\ln(2+t)^{r_3}
}{
(\frac{t}{2}+1)^{q_2}
}
\int_0^{\frac{t}{2}}
\hspace{-3mm}
\frac{
{\rm d}s
}{
s^{p_2}(1+s)^{r_2}
}
\end{equs}
and straightforward integrations.
\end{proof}
We can now complete the proof of Theorem \ref{thm:cauchy}.\par
\begin{proof}[Proof of Theorem \ref{thm:cauchy}, continued]\par
First, we recall that our goal is to prove that the map ${\cal N}$
defined by
\begin{equs}
{\cal N}[{\bf z}](t)=
\int_0^t
\hspace{-2mm}
{\rm d}s~\ed^{\L(t-s)}
\vector{0}{\partial_xh({\bf z}(s))}
\label{eqn:defNrap}
\end{equs}
satisfies $\|{\cal N}[{\bf z}]\|\leq C$ for all ${\bf z}\in{\cal
B}$ with $\|{\bf z}\|=1$. We have already proved that 
$\|\P\D{\cal N}[{\bf z}]\|_{2,\frac{3}{4}}\leq C$. The other
necessary estimates are done as follows:
\begin{equs}
\|\widehat{{\cal N}[{\bf z}]}\|_{\infty,0}&\leq
C\sup_{t\geq0}
\Bone{\frac{1}{2},\frac{1}{2}}{\frac{1}{2},\frac{1}{2}}
\leq C~,\\
\|{\cal N}[{\bf z}]\|_{2,\frac{1}{4}}&\leq
C\sup_{t\geq0}
(1+t)^{\frac{1}{4}}
\Bone{\frac{1}{2},\frac{3}{4}}{\frac{1}{2},\frac{3}{4}}\leq C~,\\
\|\P \D{\cal N}[{\bf z}]\|_{2,\frac{3}{4}}&\leq
C\sup_{t\geq0}
(1+t)^{\frac{3}{4}}
\Bone{1,\frac{3}{4}}
  {\frac{1}{2},\frac{5}{4}}
\leq C
~,\\
\|(1-\P)\D{\cal N}[{\bf z}]\|_{2,\frac{3}{4}}
&\leq \sup_{t\geq0}(1+t)^{\frac{3}{4}}
B_0[{\textstyle\frac{5}{4}}](t)\leq C~,\\
\|(1-\Q)\P\D^2{\cal N}[{\bf z}]_2\|_{2,\frac{5}{4}^{\star}}&\leq 
C\|(1-\Q)\P\D{\cal N}[{\bf z}]_2\|_{2,\frac{3}{4}}
\leq C\|\P\D{\cal N}[{\bf z}]_2\|_{2,\frac{3}{4}}\leq C~,
\label{eqn:obvious}
\\
\|\Q\P\D^2{\cal N}[{\bf z}]_2\|_{2,\frac{5}{4}^{\star}}&\leq
C
\sup_{t\geq1}{\textstyle\frac{(1+t)^{\frac{5}{4}}}{\ln(2+t)}}
~\B{\frac{3}{2},\frac{3}{4},0}
   {\frac{1}{2},\frac{5}{4},\frac{1}{2},0}
\leq C~,
\label{eqn:lastineqforremark}
\\
\|(1-\P)\D^2{\cal N}[{\bf z}]_2\|_{2,\frac{5}{4}^{\star}}
&\leq \sup_{t\geq0}(1+t)^{\frac{5}{4}}
B_0[{\textstyle\frac{5}{4}}](t)\leq C~.
\label{eqn:koversqrt}
\end{equs}
In (\ref{eqn:obvious}), we used the obvious estimates $\|\P\D
f\|_2\leq\|\P f\|_2$ and $\|(1-\Q)f\|_{2,p}\leq
2^{p-q}\|(1-\Q)f\|_{2,q}$ if $q<p$, while in
(\ref{eqn:lastineqforremark}), we made use of
${\displaystyle\sup_{|k|\leq1,t\geq0}}|k|\sqrt{1+t}\ed^{-k^2t}\leq1$,
and finally in (\ref{eqn:koversqrt}) we used $\sup_{k\in{\bf R}}
|k|(1+k^2)^{-\frac{1}{2}}=1$. Incidentally, (\ref{eqn:koversqrt})
is the only place in the above estimates where the (crucial)
presence of the extra factor $(1+k^2)^{-\frac{1}{2}}$ in the
second component of the r.h.s. of (\ref{eqn:estimateonDkernel})
is used. This concludes the proof of Theorem \ref{thm:cauchy}.
\end{proof}

\section{Remainder estimates}
\label{sect:remainderestimates}

We now make precise the sense in which the semigroup $\ed^{\L t}$
is {\em close} to that of (\ref{eqn:linearfourieruv}), whose
Fourier transform is given by
\begin{equs}
\ed^{\L_0 t}&\equiv
\myl{14}
\begin{matrix}
\ed^{-k^2t+ikt} & 0\\
0 & \ed^{-k^2t-ikt}
\end{matrix}
\myr{14}
\label{eqn:defT}~.
\end{equs}

\begin{lemma}
\label{lem:closetoheat}
Let $\P$ be the Fourier multiplier with the
characteristic function on $[-1,1]$, and let $\ed^{\L t}$
resp.~$\ed^{\L_0 t}$ be as in (\ref{eqn:defeLt}),
resp.~(\ref{eqn:defT}) and $\S$ be as in (\ref{eqn:defS}). Then
one has the estimates
\begin{equs}
\sup_{t\geq0,k\in{\bf R}}
\sqrt{1+t}\ed^{\frac{k^2t}{2}}
\my{|}{15}
\left(
\P\S
\ed^{\L t}-\ed^{\L_0 t}\S
\vbox to 13pt{}
\right)_{i,j}
\my{|}{15}
&\leq C~,
\label{eqn:closetoheat}
\end{equs}
where $(\P\S\ed^{\L t}-\ed^{\L_0 t}\S)_{i,j}$ denotes the
$(i,j)$-entry in the matrix $\P\S\ed^{\L t}-\ed^{\L_0 t}\S$.
\end{lemma}

\begin{proof}
The proof follows by considering separately $|k|\leq1$ and
$|k|>1$. We first rewrite
\begin{equs}
\P\S\ed^{\L t}-\ed^{\L_0 t}\S=
\P\myl{11}\S\ed^{\L t}-\ed^{\L_0 t}\S\myr{11}
+(1-\P)\ed^{\L_0 t}\S~.
\end{equs}
We then have
\begin{equs}
\sup_{t\geq0,k\in{\bf R}}
\sqrt{1+t}\ed^{\frac{k^2t}{2}}
\my{|}{15}
\left(
(1-\P)\ed^{\L_0 t}\S
\vbox to 13pt{}
\right)_{i,j}
\my{|}{15}
&\leq 
\sup_{t\geq0,|k|\geq1}
\sqrt{1+t}\ed^{-\frac{k^2t}{2}}
\leq C~.
\end{equs}
For $|k|\leq1$, we first compute
\begin{equs}
\ed^{\L_0 t}\S
&=
\ed^{-k^2t}
\myl{14}
\begin{matrix}
\ed^{ikt} & \ed^{ikt}\\
\ed^{-ikt} & -\ed^{-ikt}
\end{matrix}
\myr{14}~,\\
\S\ed^{\L t}
&=\ed^{-k^2t}
\myl{20}
\begin{matrix}
\cos(kt\Delta)+\frac{1-ik}{\Delta}\,i\,\sin(kt\Delta) &
\cos(kt\Delta)+\frac{1+ik}{\Delta}\,i\,\sin(kt\Delta)\\[2mm]
\cos(kt\Delta)-\frac{1+ik}{\Delta}\,i\,\sin(kt\Delta) &
-(\cos(kt\Delta)-\frac{1-ik}{\Delta}\,i\,\sin(kt\Delta))
\end{matrix}
\myr{20}~,
\end{equs}
where we recall that $\Delta=\sqrt{1-k^2}$. We next note
that
\begin{equs}
\P|\sin(kt\Delta)-\sin(kt)|+\P
|\cos(kt\Delta)-\cos(kt)|&\leq 
\P|\cos(kt(\Delta-1))-1|+\P|\sin(kt(\Delta-1))|\\
&\leq 
\P|\sqrt{1-k^2}-1|~|k|t\leq \P|k|^3t~,\\
\P|({\textstyle\frac{1}{\Delta}}-1)\sin(kt\Delta)|
&\leq
\P|\sqrt{1-k^2}-1|~|k|t\leq \P|k|^3t~.
\end{equs}
The proof is completed noting that 
\begin{equs}
\sup_{|k|\leq1,t\geq0}
t^{\frac{m}{2}}
|k|^n\ed^{-\frac{k^2t}{2}}\leq C(n)
\end{equs}
for any (finite) $0\leq m\leq n$.
\end{proof}

We are now in a position to prove that the remainder 
\begin{equs}
{\cal R}[{\bf z}](t)=
\myl{10}
\S\ed^{\L t}-\ed^{\L_0 t}\S
\myr{10}{\bf z}_0+
\int_0^t
\hspace{-2mm}
{\rm d}s~
\my{[}{14}
\S\ed^{\L(t-s)}
\vector{0}{\partial_xh({\bf z}(s))}
-
\ed^{\L_0(t-s)}\S~
\vector{0}{\partial_xg_0({\bf z}(s))}
\my{]}{14}
\end{equs}
satisfies improved estimates as stated in
(\ref{eqn:onRannounce}):

\begin{theorem}
\label{thm:simplificator}
Let $\epsilon_0$ be again the (small) constant provided by
Theorem \ref{thm:cauchy}. Then for all ${\bf z}_0\in{\cal B}_0$ with
$|{\bf z}_0|\leq\epsilon_0$, the solution ${\bf z}$ of
(\ref{eqn:p-system}) satisfies
\begin{equs}
\|{\cal R}[{\bf z}]\|_{2,\frac{3}{4}^{\star}}+
\|\D{\cal R}[{\bf z}]\|_{2,\frac{5}{4}^{\star}}\leq C\epsilon_0~.
\label{eqn:onR}
\end{equs}
\end{theorem}

\begin{proof}
We first note that
\begin{equs}
\myl{10}\S\ed^{\L t}-\ed^{\L_0 t}\S\myr{10}{\bf z}_0=
\myl{10}
\S\P\ed^{\L t}-\ed^{\L_0 t}\S
\myr{10}{\bf z}_0+\S(1-\P)\ed^{\L t}{\bf z}_0
\equiv L_1[{\bf z}_0](t)+L_2[{\bf z}_0](t)
~,
\end{equs}
and then use the fact that by Lemma \ref{lem:closetoheat}, we have
\begin{equs}
\|\D^{\alpha}L_1[{\bf z}_0]\|_{2,\frac{3}{4}+\frac{\alpha}{2}}&\leq
C\sup_{t\geq0}(1+t)^{\frac{1}{4}+\frac{\alpha}{2}}
\min\myl{10}
\|\D^{\alpha}{\bf z}_0\|_2~,~
t^{-\frac{1}{4}-\frac{\alpha}{2}}\|\widehat{{\bf z}_0}\|_{\infty}
\myr{10}
\leq C|{\bf z}_0|
\end{equs}
for $\alpha=0,1$ and finally
\begin{equs}
\|L_2[{\bf z}_0]\|_{2,\frac{3}{4}}+
\|\D L_2[{\bf z}_0]\|_{2,\frac{5}{4}}
&\leq C
(\|{\bf z}_0\|_2+\|\D{\bf z}_0\|_2)
\sup_{t\geq0}
(1+t)^{\frac{5}{4}}\ed^{-\frac{t}{4}}
\leq C|{\bf z}_0|~.
\end{equs}
This proves
\begin{equs}
\|\myl{10}\S\ed^{\L t}-\ed^{\L_0 t}\S\myr{10}{\bf z}_0\|_{2,\frac{3}{4}}+
\|\D 
\myl{10}\S\ed^{\L t}-\ed^{\L_0 t}\S\myr{10}
{\bf z}_0\|_{2,\frac{5}{4}}
\leq C|{\bf z}_0|
\end{equs}
for all ${\bf z}_0\in{\cal B}_0$. We then show that
\begin{equs}
\|
{\cal R}[{\bf z}](t)-\myl{10}\S\ed^{\L t}-\ed^{\L_0 t}\S\myr{10}{\bf z}_0
\|_{2,\frac{3}{4}^{\star}}+
\|\D\myl{11}
{\cal R}[{\bf z}](t)-\myl{10}\S\ed^{\L t}-\ed^{\L_0 t}\S\myr{10}{\bf z}_0
\myr{11}\|_{2,\frac{5}{4}^{\star}}
&\leq C\|{\bf z}\|^2
\end{equs}
for all ${\bf z}\in{\cal B}$. We only need to prove the estimates
for $\|{\bf z}\|=1$. We first
decompose
\begin{equs}
{\cal R}[{\bf z}](t)-\myl{10}\S\ed^{\L t}-\ed^{\L_0 t}\S\myr{10}{\bf z}_0=
\S{\cal N}_1[{\bf z}](t)+
\S{\cal N}_2[{\bf z}](t)+
{\cal N}_3[{\bf z}](t)~,
\label{eqn:firstdecomposition}
\end{equs} 
where
\begin{equs}
{\cal N}_1[{\bf z}](t)&=(1-\P)
\int_0^t
\hspace{-2mm}
{\rm d}s~\ed^{\L(t-s)}
\vector{0}{\partial_xh({\bf z}(s))}
~,\\
{\cal N}_2[{\bf z}](t)&=\P\int_0^t
\hspace{-2mm}
{\rm d}s~\ed^{\L(t-s)}
\vector{0}{\partial_xh({\bf z}(s))-\partial_xg_0({\bf z}(s))}~,\\
{\cal N}_3[{\bf z}](t)&=\int_0^t
\hspace{-2mm}
{\rm d}s~
\myl{10}\P\S\ed^{\L(t-s)}-\ed^{\L_0(t-s)}\S\myr{10}
\vector{0}{\partial_xg_0({\bf z}(s))}~.
\end{equs}
We then recall that $h({\bf z})$ satisfies
\begin{equs}
\|h({\bf z})\|_{2,\frac{3}{4}}+
\|\D h({\bf z})\|_{2,\frac{5}{4}}
\leq C\|{\bf z}\|^2
~,
\end{equs}
which implies
\begin{equs}
\|{\cal N}_1[{\bf z}]\|_{2,\frac{3}{4}}
&\leq C\sup_{t\geq0}
(1+t)^{\frac{3}{4}}B_0[{\textstyle\frac{3}{4}}](t)\leq C~,
~~~~~
\|\D{\cal N}_1[{\bf z}]\|_{2,\frac{5}{4}}
\leq C\sup_{t\geq0}
(1+t)^{\frac{5}{4}}B_0[{\textstyle\frac{5}{4}}](t)\leq C~.
\end{equs}
Moreover $h_0(a,b)\equiv f(a,b)\partial_xb+g(a,b)-g_0(a,b)$ satisfies
\begin{equs}
\|h_0({\bf z})\|_{1,1}+
\|\D h_0({\bf z})\|_{1,\frac{3}{2}^{\star}}
\leq C\|{\bf z}\|^2~.
\end{equs}
Here, we need to consider separately $t\in[0,1]$ and
$t\geq1$ when estimating $\|\P\D{\cal N}_2[{\bf
z}]\|_{2,\frac{5}{4}^{\star}}$. Writing again $\Q$ for the
characteristic function for $t\geq1$, we find
\begin{equs}
\|\P{\cal N}_2[{\bf z}]\|_{2,\frac{3}{4}^{\star}}
&\leq C\sup_{t\geq0}
{\textstyle\frac{(1+t)^{\frac{3}{4}}}{\ln(2+t)}}
\Bone{\frac{3}{4},1}{\frac{3}{4},1}\leq C~,
\\
\|(1-\Q)\P\D{\cal N}_2[{\bf z}]\|_{2,\frac{5}{4}^{\star}}
&\leq C\sup_{0\leq t\leq1}
(1+t)^{\frac{5}{4}}
~\Bone{\frac{3}{4},\frac{3}{2}}{\frac{3}{4},\frac{3}{2}}\leq C~,
\\
\|\Q\P\D{\cal N}_2[{\bf z}]\|_{2,\frac{5}{4}^{\star}}
&\leq C\sup_{t\geq1}
{\textstyle\frac{(1+t)^{\frac{5}{4}}}{\ln(2+t)}}
~\B{\frac{5}{4},1,0}
   {\frac{3}{4},\frac{3}{2},0,1}
\leq C~.
\end{equs}
We finally note that
\begin{equs}
\|g_0({\bf z})\|_{2,\frac{3}{4}}+
\|\D g_0({\bf z})\|_{2,\frac{5}{4}}
\leq C\|{\bf z}\|^2~,
\end{equs}
and so, using Lemma \ref{lem:closetoheat}, we find
\begin{equs}
\|{\cal N}_3[{\bf z}]\|_{2,\frac{3}{4}^{\star}}
&\leq
\sup_{t\geq0}
{\textstyle\frac{(1+t)^{\frac{3}{4}}}{\ln(2+t)}}
~\B{\frac{1}{2},\frac{3}{4},\frac{1}{2}}{\frac{1}{2},\frac{3}{4},
\frac{1}{2},0}\leq C~,
\\
\|\D{\cal N}_3[{\bf z}]\|_{2,\frac{5}{4}^{\star}}
&\leq
\sup_{t\geq0}
{\textstyle\frac{(1+t)^{\frac{5}{4}}}{\ln(2+t)}}
~\B{1,\frac{3}{4},\frac{1}{2}}
   {\frac{1}{2},\frac{5}{4},\frac{1}{2},0}\leq C~.
\end{equs}
This completes the proof.
\end{proof}

It now only remains to prove the estimates (\ref{eqn:onRtilde})
and (\ref{eqn:onRtildeLip}) on the maps $\widetilde{\cal
R}_{\{u,v\}}$, where we recall that
\begin{equs}
\widetilde{\cal R}_{u}[{\bf z},{\bf R}^{N}](t)&=
 c_{+}\E_0   [h_{1,u}+h_{3,u}](t)
-c_{-}\E_{-2}[h_{1,v}+h_{3,v}](t)
+c_{3}\E_{-1}[h_{2}  +h_{4}  ](t)
~,
\\
\widetilde{\cal R}_{v}[{\bf z},{\bf R}^{N}](t)&=
c_{-}\E_0   [h_{1,v}+h_{3,v}](t)
-c_{+}\E_{ 2}[h_{1,u}+h_{3,u}](t)
-c_{3}\E_{ 1}[h_{2}  +h_{4}  ](t)
~,
\end{equs}
with
\begin{equs}[2]
\E_{\sigma}[h](t)&=
\partial_x\int_0^{t}
\hspace{-2mm}{\rm d}s~\ed^{\partial_x^2(t-s)}~
\T^{\sigma}h(s)~&~\mbox{and}\\
h_{1,u}&=R_u^N(u+u_{\star})~,~~~h_{3,u}=u_1^2~,&~~~
h_2&=(\T R_u^N)\T^{-1}\myl{12}\frac{v+v_{\star}}{2}\myr{12}
+(\T^{-1}R_v^N)\T\myl{12}\frac{u+u_{\star}}{2}\myr{12}\\
h_{1,v}&=R_v^N(v+v_{\star})~,~~~h_{3,v}=v_1^2~,&~~~
h_4&=(\T u_{\star})(\T^{-1}v_{\star})~.
\end{equs}
Here, we will only prove that
\begin{equs}
\label{eqn:onRtilderap}
\sum_{\alpha=0}^1
\|\D^{\alpha}
\widetilde{\cal R}_{\{u,v\}}
[{\bf z},{\bf R}^{N}]\|_{2,\frac{3}{4}+\frac{\alpha}{2}-\epsilon}
&\leq
C\epsilon_0
\sum_{\alpha=0}^1
\|\D^{\alpha}{\bf R}^{N}\|_{2,\frac{3}{4}+\frac{\alpha}{2}-\epsilon}
+C~.
\end{equs}
It is then straightforward to show (\ref{eqn:onRtildeLip}), namely
that the maps $\widetilde{\cal R}_{\{u,v\}}$ are Lipschitz in
their second argument; we omit the details.

To prove (\ref{eqn:onRtilderap}), we first need estimates on ${\bf
h}_1=(h_{1,u},h_{1,v})$, $h_2$, ${\bf h}_3=(h_{3,u},h_{3,v})$ and
$h_4$. We note that ${\bf u}_0=(u_0,v_0)$ and ${\bf u}_1=(u_1,v_1)$
satisfy
\begin{equs}
\|{\bf u}_0\|_{1,0}+\|{\bf u}_1\|_{1,0}+
\|\D {\bf u}_0\|_{1,\frac{1}{2}}+\|\D {\bf u}_1\|_{1,\frac{1}{2}}
&\leq
C~,
\\
\sup_{t\geq0}
 (1+t)^{\frac{3}{2}}\myl{10}|{\bf u}_0(\pm t,t)|+|{\bf u}_1(\pm t,t)|\myr{10}
+(1+t)^{2}\myl{10}|\D {\bf u}_0(\pm t,t)|+|\D {\bf u}_1(\pm t,t)|\myr{10}
&\leq C
\end{equs}
for some constant $C$; see Proposition \ref{prop:burgers}. We thus find that 
\begin{equa} \label{eqn:hi}
\|{\bf h}_{1}\|_{1,1-\epsilon}+
\|\D {\bf h}_{1}\|_{1,\frac{3}{2}-\epsilon}
+
\|h_2\|_{1,1-\epsilon}+
\|\D h_2\|_{1,\frac{3}{2}-\epsilon}&\leq 
C\epsilon_0
\sum_{\alpha=0}^1
\|\D^{\alpha}{\bf
R}^{N}\|_{2,\frac{3}{4}+\frac{\alpha}{2}-\epsilon}~,\\
\|{\bf h}_3\|_{1,1}+\|\D {\bf h}_3\|_{1,\frac{3}{2}}
+
\|h_4\|_{1,\frac{3}{2}}+
\|\D h_4\|_{2,2}
&\leq C~.
\end{equa}
The proof of (\ref{eqn:onRtilderap}) then follows from Proposition
\ref{prop:hiinhomoheat}, which implies that
\begin{equs}
\sum_{\alpha=0}^{1}
\|\D^{\alpha}\E_{\sigma}[{\bf h}_{1}]\|_{2,\frac{3}{4}+\frac{\alpha}{2}-\epsilon}+
\|\D^{\alpha}\E_{\sigma}[h_{2}]\|_{2,\frac{3}{4}+\frac{\alpha}{2}-\epsilon}
&\leq C\epsilon_0
\sum_{\alpha=0}^1
\|\D^{\alpha}{\bf
R}^{N}\|_{2,\frac{3}{4}+\frac{\alpha}{2}-\epsilon}~,\\
\sum_{\alpha=0}^{1}
\|\D^{\alpha}\E_{\sigma}[{\bf h}_{3}]\|_{2,\frac{3}{4}+\frac{\alpha}{2}^{\star}}+
\|\D^{\alpha}\E_{\sigma}[h_{4}]\|_{2,\frac{3}{4}+\frac{\alpha}{2}^{\star}}
&\leq C
\end{equs}
for any $\sigma\in\{-2,-1,0,1,2\}$ if the estimates in (\ref{eqn:hi}) are
satisfied.

\begin{proposition}
\label{prop:hiinhomoheat}
Let $\epsilon>0$ and $\sigma\in\{-2,-1,0,1,2\}$. Then there holds
\begin{equs}
\sum_{\alpha=0}^{1}
\|\D^{\alpha}\E_{\sigma}[h_{1}]\|_{2,\frac{3}{4}+\frac{\alpha}{2}-\epsilon}
&\leq 
C
\sum_{\alpha=0}^{1}
\|\D^{\alpha}h_1\|_{1,1+\frac{\alpha}{2}-\epsilon}~,
\\
\sum_{\alpha=0}^{1}
\|\D^{\alpha}\E_{\sigma}[h_{2}]\|_{2,\frac{3}{4}+\frac{\alpha}{2}^{\star}}
&\leq C
\sum_{\alpha=0}^{1}
\|\D^{\alpha}h_2\|_{1,1+\frac{\alpha}{2}}~.
\end{equs}
\end{proposition}
\begin{proof}
Let $u_i=\E_{\sigma}[h_{i}]$. Taking the Fourier transform, we
find
\begin{equs}
\widehat{u_i}(k,t)=
ik\int_0^{t}
\hspace{-2mm}{\rm d}s~\ed^{-k^2(t-s)+i\sigma ks}\widehat{h_i}(k,s)~.
\end{equs}
We can restrict ourselves to
$\sum_{\alpha=0}^{1}
\|\D^{\alpha}h_{1}\|_{1,1+\frac{\alpha}{2}-\epsilon}=1$ and
$\sum_{\alpha=0}^{1}
\|\D^{\alpha}h_{2}\|_{1,1+\frac{\alpha}{2}}=1$. Then, it
follows that
\begin{equs}
\|\D^{\alpha}u_1\|_{2,\frac{3}{4}+\frac{\alpha}{2}-\epsilon}
&\leq C
\sup_{t\geq0}
(1+t)^{\frac{3}{4}+\frac{\alpha}{2}-\epsilon}
~\Bone{\frac{3}{4}+\frac{\alpha}{2},1-\epsilon}
   {\frac{3}{4},1+\frac{\alpha}{2}-\epsilon}
\leq C~,\\
\|\D^{\alpha}u_2\|_{2,\frac{3}{4}+\frac{\alpha}{2}^{\star}}
&\leq C
\sup_{t\geq0}
{\textstyle\frac{(1+t)^{\frac{3}{4}+\frac{\alpha}{2}}}{\ln(2+t)}}
~\Bone{\frac{3}{4}+\frac{\alpha}{2},1}
   {\frac{3}{4},1+\frac{\alpha}{2}}
\leq C
\end{equs}
for $\alpha=0,1$ as claimed.
\end{proof}


\def\Rom#1{\uppercase\expandafter{\romannumeral #1}}\def\u#1{{\accent"15
  #1}}\def\cprime{$'$} \def\cprime{$'$}

\end{document}